\documentclass[11pt]{amsart}
\usepackage{amsxtra}
\usepackage{amssymb}
\usepackage{mathrsfs} 
\usepackage{enumerate}
\usepackage{array,color}
\usepackage{hyperref}
\theoremstyle{plain}

\newcommand{\cleqn}{\setcounter{equation}{0}}
\newcommand{\clth}{\setcounter{theorem}{0}}
\newcommand {\sectionnew}[1]{\section{#1}\cleqn\clth}

\newtheorem{theorem}{Theorem}[section]
\newtheorem{lemma}[theorem]{Lemma}
\newtheorem{definition-theorem}[theorem]{Definition-Theorem}
\newtheorem{proposition}[theorem]{Proposition}
\newtheorem{corollary}[theorem]{Corollary}
\newtheorem{definition}[theorem]{Definition}
\newtheorem{example}[theorem]{Example}

\newtheorem{remark}[theorem]{Remark}
\newtheorem{conjecture}[theorem]{Conjecture}
\newtheorem{notation}[theorem]{Notation}

\newcommand \bth[1] { \begin{theorem}\label{t#1} }
\newcommand \ble[1] { \begin{lemma}\label{l#1} }

\newcommand \bpr[1] { \begin{proposition}\label{p#1} }
\newcommand \bco[1] { \begin{corollary}\label{c#1} }
\newcommand \bde[1] { \begin{definition}\label{d#1}\rm }
\newcommand \bex[1] { \begin{example}\label{e#1}\rm }
\newcommand \bre[1] { \begin{remark}\label{r#1}\rm }
\newcommand \bcj[1] { \begin{conjecture}\label{j#1}\rm }

\newcommand \bnota[1] { \begin{notation}\label{n#1}\rm }
\renewcommand {\eth} { \end{theorem} }
\newcommand {\ele} { \end{lemma} }

\newcommand {\epr} { \end{proposition} }
\newcommand {\eco} { \end{corollary} }
\newcommand {\ede} { \end{definition} }
\newcommand {\eex} { \end{example} }
\newcommand {\ere} { \end{remark} }
\newcommand {\ecj} { \end{conjecture} }

\newcommand {\enota} { \end{notation} }
\newcommand \thref[1]{Theorem \ref{t#1}}
\newcommand \leref[1]{Lemma \ref{l#1}}

\newcommand \coref[1]{Corollary \ref{c#1}}

\newcommand \deref[1]{Definition \ref{d#1}}
\newcommand \exref[1]{Example \ref{e#1}}
\newcommand \reref[1]{Remark \ref{r#1}}

\DeclareMathOperator{\sech}{sech}

\newcommand{\vocab}[1]{\textit{#1}}
\def \diff{{\backslash}}

\def \Rset {{\mathbb R}}         
\def \Cset {{\mathbb C}}

\def \Zset {{\mathbb Z}}

\def \B  {{\mathcal{B}}}               

\def \Om {\Omega}

\def \Ga {\Gamma}

\def \mt  {\mapsto}

\def \wt {\widetilde}

\DeclareMathOperator \Span { {\mathrm{Span}} }

\DeclareMathOperator \ord { {\mathrm{ord}} }

\newcommand \cord{\textnormal{cord}}
\newcommand \bw{\textnormal{bw}}
\newcommand \cbw{\textnormal{cbw}}
\newcommand{\mxx}[4]{\left(\begin{array}{cc}{#1} & {#2}\\ {#3} & {#4}\end{array}\right)}

\begin{document}
\title[Matrix discrete-continuous prolate spheroidal property $\&$ bispectrality]
{Matrix valued discrete-continuous functions with the 
prolate spheroidal property and bispectrality}
\author[W. Riley Casper]{W. Riley Casper}
\thanks{The research of M.Y.  has been supported by NSF grants DMS-2131243 and DMS--2200762,
and that of I.Z. by CONICET and by VI PPIT-US}
\address{
Department of Mathematics \\
California State University Fullerton \\
Fullerton, CA 92831\\
U.S.A.
}
\email{wcasper@fullerton.edu}

\author[F. Alberto Gr\"unbaum]{F. Alberto Gr\"unbaum}
\address{
Department of Mathematics \\
University of California, Berkeley \\
Berkeley, CA 94720 \\
U.S.A.
}
\email{grunbaum@math.berkeley.edu}

\author[Milen Yakimov]{Milen Yakimov}
\address{
Department of Mathematics \\
Northeastern University \\
Boston, MA 02115 \\
U.S.A.
}
\email{m.yakimov@northeastern.edu}

\author[Ignacio Zurri\'an]{Ignacio Zurri\'an}
\address{
Departamento de Matem\'{a}tica Aplicada II \\
Universidad de Sevilla \\
Seville, Spain
}
\email{ignacio.zurrian@fulbrightmail.org}
\date{}
\begin{abstract}
Classical prolate spheroidal functions play an important role in the study of 
time-band limiting, scaling limits of random matrices, and the distribution 
of the zeros of the Riemann zeta function. 
We establish an intrinsic relationship between discrete-continuous bispectral functions and 
the prolate spheroidal phenomenon. The former functions form a vast class, parametrized by an infinite dimensional manifold, 
and are constructed by Darboux transformations from classical bispectral functions
associated to orthogonal polynomials. 
Special cases include spherical functions. We prove that all such Darboux 
transformations which are self-adjoint in a certain sense give rise to integral 
operators possessing commuting differential operators and to discrete integral operators
possessing commuting shift operators. One particularly striking implication 
of this is the correspondence between discrete and continuous pairs of commuting operators.
Moreover, all results are proved in the 
setting of matrix valued functions, which provides further advantages for applications.
Our methods rely on the use of noncommutative matrix valued Fourier algebras 
associated to discrete-continuous bispectral functions. We produce the commuting 
differential and shift operators in a constructive way with explicit upper bounds 
on their orders and bandwidths, which is illustrated with many concrete examples.
\end{abstract}
\maketitle
\sectionnew{Introduction}
\subsection{Matrix valued discrete-continuous bispectrality}
A real, $N\times N$ matrix valued function $\Psi(k,x)$ defined on $\Zset\times V$ for some open interval $V\subseteq \Rset$ is called a \vocab{discrete-continuous bispectral function} if it simultaneously satisfies a difference equation
$$\sum_{j=-m}^m A_j(k)\Psi(k+j,x) = \Psi(k,x)\Omega(x)$$
with matrix valued coefficients, as well as a matrix valued differential equation
$$\sum_{j=0}^n\frac{\partial^j\Psi(k,x)}{\partial x^j}B_j(x) = \Lambda(k)\Psi(k,x)$$
for some real, $N\times N$ matrix valued functions $B_j(x)$ and $\Omega(x)$ defined on $V$, as well as $\Lambda(k)$ and $A_j(k)$ defined on $\Zset$.
This version of bispectrality is analogous to the idea of bispectrality introduced in \cite{DG} wherein the concept was introduced, based on observations inspired
 by computerized tomography and signal processing.
In its original form, bispectral functions were bivariate meromorphic functions $\psi(z,x)$ which were families of eigenfunctions of a differential operator in the spatial variable $x$, while simultaneously being eigenfunctions of a second differential operator in the spectral variable $z$.
The natural generalization considered in this paper, discrete-continuous bispectrality, is designed to encompass many other families of bispectral functions, including classical orthogonal polynomials, exceptional orthogonal polynomials, and matrix valued orthogonal polynomials.  As such, it constitutes a case of substantial importance with a wide variety of applications, including spectral theory, special functions, random matrices, integrable systems, and representation theory.
Note that, due to the noncommutativity of matrix multiplication, our choice of the order of products in this generalization is a significant detail: we may view $\Psi(k,x)$ as a family of eigenfunctions for a left-acting shift operator and a right-acting differential operator.
We could likewise define continuous-discrete bispectral functions by making the left variable continuous and the right variable discrete, leading to left-acting differential operators and right-acting difference operators.

Continuous-continuous bispectral functions have played a central role in recent results on commuting integral and differential operators \cite{CGYZ1,CGYZ2,CGYZ3,CY1}.
Scalar valued discrete-continuous bispectrality has been explored in the scalar case in various specific situations, starting with the work of Reach in his 1987 dissertation and the subsequent papers \cite{reach2,reach3}. Later exploration of discrete-continuous bispectral functions in the context of orthogonal polynomials can be found in \cite{grunbaum1,GH,GHH}. 
A class of scalar valued discrete-continuous bispectral functions was described in terms of an adelic flag manifold by Haine and Iliev \cite{HL}. 
The study of the matrix valued continuous-continuous bispectral functions was initiated by Zubelli in \cite{Z} and in his 1989 dissertation.
Noncommutative discrete-continuous bispectrality in the context of matrix valued orthogonal polynomials is explored in detail in \cite{CY2}
where the matrix Bochner problem was resolved.
Constructions of matrix valued discrete-continuous bispectral functions using spherical functions were given in \cite{GPT},
see also \cite{DuG,GPZ}.
\subsection{The prolate spheroidal property}
\label{1.2}
As a primary motivation for bispectrality, we have the {\em{prolate spheroidal property}}, which is the property that certain naturally appearing 
integral operators admit commuting differential operators.
This property first played a key role in the analysis of Landau, Pollak and Slepian of time-band limiting \cite{LP,S,SP} in the 60's and then in the works of Mehta \cite{Mehta} 
and Tracy-Widom \cite{TW1,TW2} on scaling limits of random matrices in the 90's. Connes, Consani, Moscovici \cite{C,CC,CM} found fundamental applications of the prolate spheroidal property 
of a completely different nature. They proved that the asymptotics of the zeros of the Riemann zeta function in two different regimes can be both modelled using 
classical prolate spheroidal functions. 

All integral operators with the prolate spheroidal property that were used in the above works have kernels of the form 
\[
K (w,z) = \int_\Ga \psi(w,x) \psi(z,x) dx
\]
for a scalar valued continuous-continuous bispectral function $\psi(z,x)$ and a contour $\Ga \subset \Rset$. This observation raises the following natural questions:
\begin{enumerate}
\item What is the matrix valued discrete-continuous analog of the classical scalar valued prolate spheroidal property?
\item Can one prove a general theorem stating that matrix valued discrete-continuous bispectral functions 
give rise to integral operators and discrete  integral operators (i.e. matrices) that 
posses the prolate spheroidal property?
\end{enumerate}

It is natural to expect that the operators obtained in the solutions of these two problems will play an equally fundamental role as the one played by their continuous-continuous counterparts 
in time-band limiting, random matrices and in connection with the Riemann zeta function. These questions were raised in \cite{GPZ} where some examples of
commuting differential operators were constructed. There are two key components of this setting: 
\begin{enumerate}
\item[(i)] considering (noncommutative) matrix valued functions as opposed to scalar valued ones and 
\item[(ii)] transitioning from the continuous-continuous setting to the much less well understood one of discrete-continuous ones, where 
we have a close relation to orthogonal polynomials.
\end{enumerate}
We believe that both new features will be useful in a wide range of concrete applications.

The goal of this paper is to address both problems (1) and (2). The upshot of our results is the construction of infinite-dimensional families of operators 
satisfying the matrix valued discrete-continuous analog of the prolate spheroidal property. 

Compared to \cite{CGYZ1,CGYZ2,CGYZ3,CY1}, we overcome a number of difficulties in the present paper:
\begin{itemize}
\item All functions considered here are matrix valued, and as a consequence, all algebraic arguments need to deal with 
a noncommutative setting instead of the simpler commutative one.
\item Since one of the arguments is discrete and ranges over a subset of  $\Zset$, we can no longer resort 
on functions $\Psi(x,k)$ that are analytic on a subdomain of $\Cset^2$. 
\end{itemize}

Firstly, to each matrix valued discrete-continuous bispectral function $\Psi$ we associate the integral operator
\begin{equation}\label{eqn:cont integral}
T_\Psi: F(y)\mapsto \int_{x_0}^{x_1} F(x)K(x,y) dx,\ \ K(x,y) = \sum_{k \in I} \Psi(k, x)^*\Psi(k,y)
\end{equation}
where here $I\subseteq\Zset$ is a finite set, $(x_0,x_1)\subseteq V$, and $*$ is the matrix transposition.
Under natural assumptions on $I$ and $(x_0,x_1)$, the operator $T_\Psi$ defines an integral operator on a Hilbert space with a basis consisting of smooth functions on $[x_0,x_1]$.
In contrast to the continuous-continuous situation, we can also form a discrete integral operator
\begin{equation}\label{eqn:disc integral}
S_\Psi: F(m)\mapsto \sum_{k\in I} J(m,k)F(k),\ \ J(m,n) = \int_{x_0}^{x_1} \Psi(m, y) \Psi(k,y)^* dy.
\end{equation}
Alternatively, $S_\Psi$ may be expressed as an $|I|\times|I|$ matrix operator whose entries are $ N\times N$ matrices.
Our main theorem will establish that in wide generality the operators $T_\Psi$ and $S_\Psi$ 
commute with matrix valued differential and shift operators, respectively.
Again, the non-commutativity of matrix multiplication makes the order of the products in the above expressions an important aspect of the construction.

The prototypical example of this is the projection operator defined by the Christoffel-Darboux kernel of a sequence of orthogonal matrix polynomials $P_n(x)$ for a weight matrix $W(x)$ supported on a real interval $V$ and given by
$$T(f) = \int_{x_0}^{x_1} f(x)K(x,y)dy,\ \ K(x,y) := \sum_{k=0}^\ell U(x)^*P_k(x)^* (H_kH_k^*)^{-1}P_k(y)U(y),$$
where 
\[
H_kH_k^* = \int_V P_k(x)W(x)P_k(x)^*dx, \quad W(x) = U(x)U(x)^*
\]
are the respective Cholesky decompositions and $(x_0,x_1)\subseteq V$.
This is exactly the integral operator defined by \eqref{eqn:cont integral} for the discrete-continuous bispectral function $\Psi(k,x) = H_k^{-1}P_k(x)U(x)$, with $I=\{0,1,\dots,\ell\}$.
\subsection{Results}
Our main tool to build a bridge between matrix valued discrete-continuous bispectral functions and the prolate spheroidal property is the 
notion of Fourier algebras. The {\em{left and right Fourier algebras}} of such a bispectral function 
$\Psi(k,x)$ are the algebras consisting of $M_N(\Rset)$-valued shift operators $P(k, \mathscr S_k)$ and differential operators 
$D(x, \partial_x)$, respectively, such that 
\[
P (k, \mathscr S_k) \cdot \Psi(k,x) = \Psi(k, x) \cdot D(x, \partial_x).
\]
Here, the shift operator acts by $\big(\sum_{n=-\ell}^\ell A_n(k)\mathscr S_k^n \big) F(k) := \sum_{n=-\ell}^\ell A_n(k)F(k+n)$ for arbitrary 
functions $A_n : \Zset \to M_N(\Rset)$; the differential operator acts by the right action \eqref{r-act}.
The right action is needed in order for the two actions to commute with each other, so, algebraically, we deal with 
a bimodule for two algebras consisting of shift operators and differential operators, respectively.

When the function $\Psi(k,x)$ has trivial left and right annihilators
(which holds in broad generality), the assignment $P (k, \mathscr S_k) \mt D(x, \partial_x)$ defines an 
isomorphism between the left and right Fourier algebras of $\Psi)$, 
which will be called the {\em{generalized Fourier map}}. 
The two Fourier algebras and the isomorphism between them is a far reaching generalization of the Fourier 
transform adapted to the setting of arbitrary bispectral functions; the Fourier transform arises from the 
simplest scalar valued continuous-continuous bispectral function $\exp(xz)$.

The classical orthogonal polynomials of Hermite, Laguerre, and Jacobi give rise to instances of discrete-continuous bispectral functions, which we refer to as
{\em{classical discrete-continuous bispectral functions}}; they are eigenfunctions of a second order differential operator and a shift operator of bandwidth two, and are defined in the table in Figure \ref{classical discrete-continuous table}.
For an arbitrary positive integer $N$, we consider the matrix valued analogs of these functions by multiplying by the 
identity matrix $I_N$. A {\em{Darboux transformations}} of such a function $\Psi(k,x)$ is given by
\[
\Psi(k,x) \mapsto  \wt\Psi(k,x) = \big( F(k)^{-1} P (k, \mathscr S_k) \big) \cdot \Psi(k,x) G(x)^{-1}
\]
for some matrix valued shift operator $P$ in $k$ and matrix valued functions $F(k)$ on $\Zset$ and $\wt G(x)$ on $V$, respectively, with the property
that they define a factorization of a matrix valued shift operator for which $\Psi(k,x)$ is a generalized eigenfunction:
\[
\big( \wt{P} (k, \mathscr S_k) \wt{F}(k)^{-1} F(k)^{-1} P (k, \mathscr S_k) \big) \cdot \Psi(k,x)
= \Psi(k,x)\wt G(x){G}(x)
\]
for some matrix valued shift operator $\wt{P} (k, \mathscr S_k)$ and matrix valued functions $\wt F(k)$ on $\Zset$ and $\wt G(x)$ on $V$, respectively.
For more details, see Section 2.

Of particular importance are the {\em{bispectral}} Darboux transformations: those satisfying the property that the shift operators 
$\wt{P} (k, \mathscr S_k)$ and $P (k, \mathscr S_k)$ and the functions $F(k)\wt F(k)$, $G(x)\wt G(x)$ 
are in the left and right Fourier algebras of $\Psi(k,x)$, respectively.
By the main theorem of \cite{GHY}, all bispectral Darboux transformations are also matrix valued discrete-continuous bispectral functions.

These Darboux transformations present a method for constructing bispectral functions which vastly generalizes the method of 
Duistermaat and Gr\"{u}nbaum \cite{DG} by recursive factorizations of second order differential operators and switching their orders, 
and Wilson's algebro-geometric method based on involutions of the adelic Grassmannian \cite{Wilson}; we refer the reader to Sect. \ref{bDtr} 
for a detailed discussion. Here, we note that all known bispectral functions 
are bispectral Darboux transformations in the above sense from basic bispectral functions.

The assumption that $P(k, \mathscr S_k)$ and $\wt{P}(k, \mathscr S_k)$ 
are in the left Fourier algebra means that there exist right-acting matrix valued differential 
operators $D(x, \partial_x)$ and $\wt{D}(x, \partial_x)$ satisfying
\[
P (k, \mathscr S_k) \cdot\Psi(k,x) = \Psi(k,x)\cdot D(x,\partial_x)
\quad\text{and}\quad \wt{P} (k, \mathscr S_k) \cdot\Psi(k,x) = \Psi(k,x)\cdot \wt D(x,\partial_x).
\]
The differential operators $D$ and $\wt{D}$ are the images of the shift operators $P$ and $\wt{ P}$ 
under the generalized Fourier isomorphism of the bispectral function $\Psi(k,x)$.
We call a bispectral Darboux transformation $\Psi \mt {\wt\Psi}$ {\em{self-adjoint}} when
\[ 
P(k, \mathscr S_k)^* = \wt{P} (k, \mathscr S_k), \; \; D(x,\partial_x)^* = \wt D(x,\partial_x) \quad \mbox{and} \quad
F(k)^* = \wt F(k), \; \; G(x)^* = \wt G(x),
\]
where star denotes the formal adjoint of matrix valued differential and shift operators (extending transposition of matrices).

We define the {\em{degree}} of the Darboux transformation ${\wt\Psi}$ to be the pair $(d_1, d_2)$ where $d_1$ is the bandwidth of the shift operator 
$P(k, \mathscr S_k)$ and 
$d_2$ is the order of the image of $P (k, \mathscr S_k)$ under the generalized Fourier map (which is a matrix valued differential operator in $x$). 

Our main theorem resolves problems (1) and (2) in Section \ref{1.2} and provides an effective upper 
bound on the orders and bandwidths of commuting differential and shift operators:
\smallskip

\noindent
{\bf{Main Theorem.}} {\em{Let ${\wt\Psi}(k,x)$ be a self-adjoint bispectral Darboux transformation of 
a classical discrete-continuous bispectral function $\Psi(k,x)$ of degree $(d_1,d_2)$ 
with the property that the leading terms of the operators $P (k, \mathscr S_k)$ and $D(x,\partial_x)$ are nondegenerate.
Then for $I=\{0,1,\dots, n\}$ or $I=\{-1,-2,\dots,-n\}$ and $(x_0,x_1)$ a suitably chosen interval in the domain of $\Psi$ the following hold:
\begin{enumerate}
\item The continuous integral operator $T_{\wt\Psi}$ defined by \eqref{eqn:cont integral} 
commutes with a matrix valued differential operator of positive order $\leq 2 d_1 d_2$
which belongs to the right Fourier algebra of ${\wt\Psi}$.  
\item The discrete integral operator $S_{\wt\Psi}$ defined by \eqref{eqn:disc integral} 
commutes with a matrix valued shift operator of positive order $\leq 2 d_1 d_2$
which belongs to the left Fourier algebra of ${\wt\Psi}$. This shift operator is the inverse image under 
the generalized Fourier map of the differential operator in part (1) of the theorem.
\end{enumerate}
For the bounds, we take:
$x_0=-\infty$ or $x_1=\infty$ if $\Psi$ is of Hermite type;
$x_0=0$ or $x_1=\infty$ if $\Psi$ is of Laguerre type;
and either $x_0=-1$ or $x_1=1$ if $\Psi$ is of Jacobi type.
}}
\medskip

\noindent
We note the following: 
\medskip

(1) The actual theorem proved in Section \ref{6} is stronger and deals with a {\em{robustness assumption}} that relaxes the assumption that 
the operators $P(k, \mathscr S_k)$ and $D(x,\partial_x)$ defining the Darboux transformation 
$\Psi(k,x) \mapsto {\wt\Psi}(k,x)$ have  
nondegenerate leading terms. The robustness is only a {\em{genericity assumption}} and is always satisfied if the size $N$ equals $1$.
\medskip

(2) Although the original bispectral functions $\Psi$ are scalar multiples of the identity matrix $I_N$, the Darboux transformations ${\wt\Psi}(k,x)$ 
are very general matrix valued functions that are not in any way scalar multiples of the identity matrix (except in the cases of
very special Darboux transformations). We refer the reader to Section \ref{7.3} for an example about how Darboux transformations 
produce complicated (full matrix) bispectral functions. 
\medskip

\noindent
{\bf{Remark.}} By using the methods of \cite[Theorems 6.8 and 6.10]{CY1}, one can prove that the self-adjoint bispectral Darboux transformations
of classical discrete-continuous bispectral functions are parametrized by the points of infinite dimensional Grassmannians
generalizing Wilson's adelic Grassmannain \cite{Wilson}.
\medskip

(3) The theorem provides an effective way of constructing commuting differential and shift operators for the integral and discrete integral operators in question. 
This is done in two steps. The first one is describing the corresponding right and left Fourier algebras of the discrete-continuous bispectral function 
${\wt\Psi}(k,x)$ of degrees/bandwidths $\leq 2 d_1 d_2$. The second is solving linear equations coming from bilinear concomitants. 
\medskip

\noindent
{\bf{Correspondence.}} A striking feature of the main theorem is that it creates a correspondence between 
pairs of discrete commuting operators \eqref{eqn:disc integral} on the one side, and pairs of continuous commuting operators \eqref{eqn:cont integral} on the other:
\[
\big( T_{\wt\Psi}, R(x,\partial_x) \big)  \longleftrightarrow \big(S_{\wt\Psi}, (b_{\wt\Psi}^{-1} R) (k, \mathscr S_k) \big).
\]
These pairs are parametrized by the self-adjoint bispectral Darboux transformations $\wt\Psi(k,x)$ of classical discrete-continuous bispectral functions 
(forming infinite dimensional manifolds). On the left side we have a commuting pair of an integral operator and an operator in the right Fourier algebra of  $\wt\Psi(k,x)$. 
On the right side 
we have a discrete integral operator and the inverse image of the differential operator in the first pair under the generalized Fourier map 
associated to  $\wt\Psi(k,x)$.
\medskip

The utility of the main theorem is demonstrated in Section \ref{7}, where we first provide a simple and brief derivation of the main results of \cite{grunbaum1}
based on intrinsic arguments with Fourier algebras.
Following this, we provide a brand new example of a commuting family of integral and differential operators associated to a self-adjoint 
bispectral Darboux transformation of the Laguerre polynomials.  Next, we obtain a matrix valued example, associated to a noncommutative bispectral 
Darboux transformation of the Hermite polynomials. Lastly, we conclude with an example of discrete-continuous bispectral functions which are derived 
from soliton solutions of the KdV equation, instead of coming from classical orthogonal polynomials.  
We show again that the the associated integral operator commutes with a differential operator. One can produce many additional matrix valued examples 
but their complexity grows rapidly and they cannot be easily displayed on paper.

The paper is organized as follows. Section \ref{2} contain background material on matrix valued discrete-continuous bispectrality and Fourier algebras attached to bispectral functions. 
Section \ref{3} describes the Hermite, Laguerre, and Jacobi orthogonal polynomial which we use and their treatment in the 
setting of discrete-continuous bispectral functions. In Section \ref{4} we establish sharp upper bounds on the sizes of certain canonical bifiltrations of the Fourier algebras 
of all bispectral Darboux transformations from classical discrete-continuous bispectral functions. Section \ref{5} proves results on the bilinear concomitants of matrix valued differential 
and shift operators that are in turn used to construct commuting differential and shift operators for integral and discrete integral operators, respectively. 
The main theorem is proved in Section \ref{6}. 
Section \ref{7} contains examples illustrating the broad applications of the main theorem.
\medskip

{\bf{Notation and conventions.}} All vector spaces in the paper will be over the real numbers and $\dim$ will denote $\dim_\Rset$.
By a smooth function we mean a $C^\infty$-function.

\sectionnew{Bispectrality and Fourier algebras}
\label{2}
The classical bispectral functions explored by Duistermaat and Gr\"{u}nbaum \cite{DG} were bivariate functions $\psi(k,x)$ which were simultaneously families of eigenfunctions for differential operators in each variable, i.e. 
\begin{align}
L(x,\partial_x)\cdot \psi(k,x) &= \lambda(k)\psi(k,x),
\label{bisp1}
\\
B(k,\partial_k)\cdot \psi(k,x) &= \theta(x)\psi(k,x).
\label{bisp2}
\end{align}
The associated differential operators $L(x,\partial_x)$ and $B(k,\partial_k)$ are called \vocab{bispectral operators}.
As a key insight toward classifying all rank one bispectral fuctions, Wilson initiated the study of the algebras of \emph{all} bispectral operators associated to $\psi(k,x)$,
which form the so called \vocab{bispectral algebras} \cite{Wilson}.
As a further insight, \cite{BHY1,GHY} considered the \vocab{Fourier algebras} consisting of differential operators $P(k,\partial_k)$ and $B(x,\partial_x)$ satisfying 
\[
P(k,\partial_k)\cdot\psi(k,x) = D(x,\partial_x)\cdot\psi(k,x),
\]
using them to study Darboux transformations in greater generality \cite{GHY}.
Further study of bispectral and Fourier algebras has led to many recent important advances in the study of bispectrality, prolate spheroidal operators and special functions \cite{CGYZ1,CGYZ2,CY1,CY2}.

This section provides a general abstract approach to bispectrality via bimodules for noncommutative algebras. Each such function 
gives rise to an associated pair of Fourier algebras and a canonical isomorphism between them. 
Fourier algebras are in turn at the heart of a construction of  
new bispectral functions from old ones. 
While we will provide definitions in complete generality, we will soon specify to matrix valued differential and shift operators.
However, many of the results in this paper will apply to more general operator algebras, such as differential or shift operators taking values in a more general finite-dimensional real algebra.
\subsection{Basic definitions}
Bispectrality arises in the context of operator algebras, such as algebras of differential or shift operators, acting on certain function spaces.
For our purposes, we adopt the following definition of an operator algebra.
\bde{operator algebra}
By an \vocab{operator algebra}, we mean a real algebra $\mathcal{D}$ with an $\Rset$-linear anti-involution $*$ and with a distinguished subalgebra 
$\mathcal{K}$ of \vocab{multiplicative operators}. The operation $*$ is referred to as the \vocab{formal adjoint} of $\mathcal D$.
\ede
This differs slightly from the usual definition of an operator algebra, in that we are avoiding explicit mention of a topological vector space on which our operators act linearly.
We call the operators in $\mathcal{K}$ multiplicative because in all important examples $\mathcal{K}$ will consist of multiplication operators by functions, which will 
be used for writing spectral equations. 
Following \cite[\S 2]{BHY1} and \cite[\S 2.1]{GHY}, we define the abstract notion of a bispectral function: 
\bde{bisp context}
A \vocab{bispectral context} is a triple $(\mathcal D_k,\mathcal D_x,\mathcal M)$ with $\mathcal D_k$ and $\mathcal D_x$ operator algebras and $\mathcal M$ a 
$(\mathcal D_k,\mathcal D_x)$-bimodule.
A bispectral context is called \vocab{commutative} if the associated subalgebras $\mathcal K_k,\mathcal K_x$ of multiplicative operators are commutative, and is called \vocab{noncommutative} otherwise.
A \vocab{bispectral triple} for this bispectral context is a triple $(B_k,L_x,\Psi)$ with $B_k\in \mathcal D_k\diff\Rset$, $L_x \in \mathcal D_x \diff\Rset$ 
and $\Psi\in \mathcal M$ with trivial left and right annihilators satisfying
$$B_k\cdot \Psi = \Psi\cdot G_x\ \ \text{and}\ \ \Psi\cdot L_x = F_k\cdot\Psi$$
for some $F_k\in \mathcal K_k\diff \Rset$ and $G_x\in \mathcal K_x\diff\Rset$.
In this case $\Psi\in \mathcal M$ is called a \vocab{bispectral function}.
A bispectral triple $(B_k, L_x,\Psi)$ is called \vocab{self-adjoint} if $B_k = B_k^*$, $L_x=L_x^*$, $G_x^*=G_x$, and $F_k^*=F_k$, in which case the bispectral function $\Psi$ is also called self-adjoint.
\ede

\bde{Fourier alg}
Suppose that $(\mathcal D_k,\mathcal D_x,\mathcal M)$ is a bispectral context and $\Psi\in \mathcal M$ is a bispectral function.
The \vocab{left and right Fourier algebras} of $\Psi$ are the $\Rset$-subalgebras of $\mathcal D_k$ and $\mathcal D_x$ defined by
\begin{align*}
\mathcal F_k(\Psi) &:= \{P_k\in\mathcal D_k : P_k\cdot\Psi = \Psi\cdot D_x\ \text{for some $D_x\in\mathcal D_x$}\}, 
\\
\mathcal F_x(\Psi) &:= \{D_x\in\mathcal D_x : P_k\cdot\Psi = \Psi\cdot D_x\ \text{for some $D_k\in\mathcal D_k$}\}.
\end{align*}
These subalgebras are further refined into the \vocab{left and right bispectral algebras} of $\Psi$, defined by
\begin{align*}
\mathcal B_k(\Psi) &:= \{P_k\in\mathcal D_k : P_k\cdot\Psi = \Psi\cdot D_x\ \text{for some $D_x\in\mathcal K_x$}\},
\\
\mathcal B_x(\Psi) &:= \{D_x\in\mathcal D_x : P_k\cdot\Psi = \Psi\cdot D_x\ \text{for some $D_k\in\mathcal K_k$}\}.
\end{align*}
\ede
The left and right annihilators of $\Psi$ are trivial by definition, and as a consequence the left and right Fourier algebras of $\Psi$ are isomorphic via
$$b_\Psi: \mathcal F_k(\Psi)\rightarrow \mathcal F_x(\Psi),\ \ b_\Psi(P_k) := D_x\ \ \text{with}\ \ P_k\cdot\Psi = \Psi\cdot D_x.$$
\bde{Fourier map}
The algebra isomorphism $b_\Psi$ is called the \vocab{generalized Fourier map}.
\ede

We do no assume any compatibility conditions on $b_\Psi$ with respect to the $*$ anti-involutions of 
$\mathcal D_x$ and $\mathcal D_k$. In an important situation such compatibility holds, 
see \reref{b-star} below.

\bre{Fourier}
The name generalized Fourier map comes from the simplest classical case of bispectrality, when $\mathcal D_k$ is the algebra of differential operators in variable $k$, $\mathcal D_x$ is the opposite algebra of differential operators in variable $x$ and $\Psi = e^{kx}$.
Then $\mathcal F_k = \Rset[k,\partial_k]$ and $\mathcal F_x=\Rset[x,\partial_x]^{op}$, where the latter is the opposite algebra to the first Weyl 
algebra $\Rset[x,\partial_x]$, i.e. the algebra generated by $\partial_x$ and $x$ subject to the relation $\partial_x x - x \partial_x = -1$.
The isomorphism $b_\Psi$ is the Fourier transform which sends $k\mapsto \partial_x,\ \partial_k\mapsto x$.
\ere

In this paper, a key role will be played by the subspaces of $\mathcal F_k(\Psi)$ and $\mathcal F_x(\Psi)$ consisting of {\em{formally bisymmetric operators}}, 
defined by $\mathcal F_{x,sym}(\Psi) := b_\Psi(\mathcal F_{k,sym}(\Psi))$ and
$$\mathcal F_{k,sym}(\Psi) := \{ P_k \in\mathcal F_k(\Psi): P_k^* = P_k\ \text{and}\ b_\Psi(P_k)^*=b_\Psi(P_k)\}.$$
Note that they are $\Rset$-subspaces, but generally not $\Rset$-subalgebras, because $*$ is an anti-involution of 
$\mathcal D_x$ and $\mathcal D_k$.

\subsection{Bispectral Darboux transformations}
\label{bDtr}
New bispectral functions can be obtained from old ones by means of bispectral Darboux transformations.
This procedure was first used in \cite{DG} were scalar valued second order differential operators were obtained by recursive 
Darboux transformations from the Bessel operators
\[
L_0(x, \partial_x) := \partial_x^2 + \nu (\nu +1)/x^2 \curvearrowright L_1(x, \partial_x) \curvearrowright \ldots \curvearrowright
L_n(x, \partial_x).
\]
The Darboux process amounted to a factorization into first order differential operators and then interchanging their order
\[
L_j (x, \partial_x) = P_j(x, \partial_x) Q_j(x, \partial_x) \curvearrowright 
L_{j+1} (x, \partial_x) := Q_j(x, \partial_x) P_j(x, \partial_x). 
\]
On the level of bispectral functions, the eigenfunctions of these operators are related by
\[
\psi_j (x,k) \curvearrowright \psi_{j+1}(x,k) := P_j(x, \partial_x) \psi_j(x,k).
\] 
There were two key problems here:
\begin{enumerate}
\item It is a highly nontrivial fact that each $\psi_j(x,k)$ is an eigenfunction in the variable $k$ as in \eqref{bisp2}.
\item In order for this to hold, the Darboux transformations require for the operators 
$P_j(x, \partial_x)$ and $Q_j(x, \partial_x)$ to have rational coefficients 
in the KdV case $\nu \in \Zset$ and to have rational coefficients and to be invariant under the transformation 
$x \mt - x$ in the even case $\nu \in 1/2 + \Zset$. 
\end{enumerate}

The Darboux transformation process, as a tool of constructing bispectral functions, took on a new and much more general form
in the works \cite{BHY2,KR} where all operators of the spectral algebra of a bispectral function $\Psi_0(x,k)$ 
were used for factorization. This allowed all Darboux transformations to be obtained in one single step 
and much more general bispectral functions to be obtained in this fashion. Furthermore, the algebro-geometric
method of Wilson \cite{Wilson}, based on involutions of the adelic Grassmannian to handle problem (1) above 
(in the much more general case of rank 1 continuous-continuous bispectral functions) was fully phrased 
in terms of Darboux transfomations from the whole spectral algebra.

Eventually, problems (1) and (2) 
were fully resolved in \cite{BHY1,GHY} through the method of factorizations in localizations of Fourier algebras. 
The different nature of the factorizations in problem (2) are simply the shadow of the different nature of the Fourier 
algebras in the two cases (KdV vs even case) and the dual bispectral equation in problem (1) is the result of the 
application of the Fourier isomorphism between the left and right Fourier algebras. The noncommutative algebra 
arguments needed to pass from the scalar valued situation to the matrix valued one were carried out in \cite{GHY}.

We next review the fundamentals of this method:

\bde{bispectral}
Let $\Psi$ be a bispectral function in a bispectral context $(\mathcal D_k,\mathcal D_x,\mathcal M)$, i.e., $\mathcal D_k$ and $\mathcal D_x$ operator algebras and $\mathcal M$ a 
$(\mathcal D_k,\mathcal D_x)$-bimodule.
\begin{enumerate}
\item A \vocab{bispectral Darboux transformation} $\wt\Psi$ of $\Psi$ is an element of $\mathcal M$ satisfying
$$\wt \Psi = F^{-1} P_k\cdot \Psi G^{-1}\ \ \text{and}\ \ \Psi = \wt P_k\wt F^{-1}\cdot\wt\Psi\wt G^{-1}$$
for some $P_k,\wt P_k\in \mathcal F_k(\Psi)\subset \mathcal D_k$, and some multiplicative operators $F,\wt F\in \mathcal K_k^\times$ and $G,\wt G\in\mathcal K_x^\times$
with $D_xG^{-1}\wt G^{-1}\wt D_x \in \mathcal B_x(\Psi)$ and that $\wt P_k$ is not a left zero divisor of $\mathcal D_k$.

\item A bispectral Darboux transformation is called \vocab{self-adjoint} if $F(k)^* = \wt F(k)$, $G(x)^* = \wt G(x)$, 
$P_k^* = \wt P_k$, and $b_\Psi(P_k)^* = b_\Psi(\wt P_k)$.
\end{enumerate}
\ede

Here and below, for a (noncommutative) $\Rset$-algebra $\mathcal K$, we denote by $\mathcal K^\times$ its multiplicative group of invertible elements.
We have:
\bth{bisp}
Let $\wt\Psi$ be a bispectral Darboux transformation of a bispectral function $\Psi$.
Then: 
\begin{enumerate}
\item $\wt\Psi$ is a bispectral and, more precisely, 
\begin{equation*}\label{darboux relations}
\left\lbrace\begin{array}{cc}
\wt P_k\wt F^{-1}F^{-1} P_k\cdot\Psi &= \Psi\wt GG\\ \nonumber
\Psi\cdot D_x G^{-1}\wt G^{-1}\wt D_x &= F \wt F \Psi\\ \nonumber
F^{-1} P_k\wt P_k\wt F^{-1}\cdot \wt\Psi &= \wt\Psi G\wt G\\
\wt\Psi\cdot \wt G^{-1}\wt D_x D_x G^{-1} &= \wt F F \wt\Psi \nonumber
\end{array}\right.
\end{equation*}
where $D_x := b_\Psi(P_k)$ and $\wt{D}_x := b_\Psi(\wt P_k)$.
\item If the bispectral Darboux transformation is self-adjoint, then both $\Psi$ and $\wt\Psi$ are self-adjoint.
\end{enumerate}
\eth
\begin{proof}
The first and third identities in the statement of the theorem follow easily from the definition of a bispectral Darboux transformation.
The fourth identity is easily derived from the second and the statements about self-adjointness follow immediately.
The above identities also guarantee that $\wt\Psi$ has trivial left and right annihilators.

To prove the second identity, we start with first identity and multiply both sides on the right by $G^{-1}\wt G^{-1}\wt D_x$.  Then we use the fact that both $P_k\cdot\Psi = \Psi\cdot D_{x}$ and $\wt P_k\cdot\Psi = \Psi\cdot \wt D_{x}$ to obtain
$$(\wt P_k\wt F^{-1}F^{-1}\cdot \Psi)\cdot D_x G^{-1}\wt G^{-1}\wt D_x = \wt P_k\Psi.$$
Therefore
$$(\wt P_k\wt F^{-1}F^{-1})(b_\Psi^{-1}(D_xG^{-1}\wt G^{-1}\wt D_x) - F\wt F)\cdot \Psi = 0.$$ 
Since $\Psi$ has trivial annihilators, this gives the identity,
$$(\wt P_k\wt F^{-1}F^{-1})(b_\Psi^{-1}(D_xG^{-1}\wt G^{-1}\wt D_x) - F\wt F) = 0.$$
Then since $\wt P_k$ is not a left zero divisor, this proves our second identity.
\end{proof}

\bre{dual} The notion of a bispectral Darboux transformations is self dual in the sense that it can be phrased 
both in terms of factorizations in the algebra $\mathcal D_k$ and $\mathcal D_x$. More precisely, 
\deref{bispectral} is equivalent to saying that 
\begin{equation}
\label{2nd-set-indentities}
\wt \Psi = F^{-1} \Psi \cdot D_x G^{-1}\ \ \text{and}\ \ \Psi =  \wt F^{-1} \wt\Psi \cdot  \wt G^{-1} \wt D_x
\end{equation}
for some $D_k,\wt D_x \in \mathcal F_x(\Psi)$, $F,\wt F\in \mathcal K_k^\times$ and $G,\wt G\in\mathcal K_x^\times$.

Due to the symmetric nature of these identities vs those in \deref{bispectral}, we only have to show that the identities in 
\deref{bispectral} imply \eqref{2nd-set-indentities}. In the setting of \deref{bispectral}, denote 
$D_x := b_\Psi(P_k)$ and $\wt{D}_x := b_\Psi(\wt P_k)$. Then 
\[
\wt \Psi = F^{-1} P_k \cdot \Psi G^{-1} = F^{-1} \Psi \cdot D_x G^{-1}.
\]
Using this identity and then the second identity in \thref{bisp}(1) gives
\[
\wt F^{-1} \wt\Psi \cdot \wt G^{-1} \wt D_x =
\wt F^{-1}  F^{-1} \Psi \cdot D_x G^{-1}  \wt G^{-1} \wt D_x = \wt F^{-1}  F^{-1} F \wt F \Psi = \Psi,
\]
which proves that that the identities in \deref{bispectral} are equivalent to those in \eqref{2nd-set-indentities}.
\ere

\subsection{Shift and differential operators}
In this paper, we will be concerned only with the single bispectral context $(\Sigma(\Zset),\Omega(V), M_N(C^\infty(\Zset\times V)))$, 
which we now define.

\bde{shift operators}
A \vocab{difference operator} or \vocab{shift operator} is an operator of the form 
$L (k, \mathscr S_k) := \sum_{n=-\ell}^\ell A_n(k)\mathscr S_k^n$ for some functions $A_n: \mathbb Z\rightarrow M_N(\Rset)$.
It acts on $M_N(\Rset)$-valued functions $F: \mathbb{Z} \rightarrow M_N(\Rset)$ by
$$L (k, \mathscr S_k) \cdot F(k) := \sum_{n=-\ell}^\ell A_n(k)F(k+n).$$
\ede
In other words, $\mathscr S_k^j$ is the basic \vocab{$k$-shift operator}, acting on matrix valued functions $F: \mathbb{Z} \rightarrow M_N(\Rset)$ 
by $\mathscr S_k^j\cdot F(k) := F(k+j)$.

The collection of all shift operators forms an $\Rset$-algebra $\Sigma(\Zset)$ with a product satisfying the fundamental relation
$\mathscr S_k^n A(k) = A(k+n)\mathscr S_k^n$ for all $n \in \Zset$ and functions $A: \mathbb{Z}\rightarrow M_N(\Rset)$.
Consider the anti-involution $*$ on $\Sigma(\Zset)$ given by
$$\left(\sum_{n=-\ell}^\ell A_n(k)\mathscr S_k^n \right)^* := \sum_{n=-\ell}^\ell A_n(k-n)^*\mathscr S_k^{-n},$$
where in the right hand side $*$ denotes the matrix transpose.
This anti-involution makes $\Sigma(\Zset)$ an operator algebra whose subalgebra of multiplicative operators consists 
of matrix valued functions on $\mathbb Z$.

Next, let $V\subseteq \Rset$ be an open interval and $\Omega(V)$ be the opposite algebra 
of the algebra of differential operators (see \reref{Fourier}) with $M_N(\Rset)$-valued smooth coefficients on $V$. 
It has a canonical right action on the algebra of $M_N(\Rset)$-valued smooth functions on $V$. An operator of the form 
$D(x,\partial_x) = \sum_{n=0}^\ell \partial_x^n B_n (x)$ acts on a smooth function $F: V\rightarrow M_N(\Rset)$ by
\begin{equation}
\label{r-act}
F\cdot D(x,\partial_x)  := \sum_{n=0}^\ell \frac{d^n F(x)}{dx^n }B_n (x).
\end{equation}
The anti-involution $*$ of $M_N(\Rset)$ extends to an anti-involution of $\Omega(V)$ defined by
$$
\left(\sum_{n=0}^m \partial_x^n B_n (x)\right)^* := \sum_{n =0}^m (-1)^n B_n (x)^*\partial_x^n
$$
for $*$ the matrix transpose in the right hand side.
This makes $\Omega(V)$ an operator algebra whose subalgebra of multiplicative operators 
consists of matrix valued smooth functions on $V$.

The collection $M_N(C^\infty (\mathbb{Z}\times V))$ of $M_N(\Rset)$-valued functions on $\mathbb{Z}\times V$ smooth in the second variable 
is a $(\Sigma(\Zset),\Omega(V))$-bimodule 
with the natural action defined above, and $(\Sigma(\Zset),\Omega(V), M_N(C^\infty(\Zset \times V)))$ forms a bispectral context.
\bde{dc context}
We call the bispectral context $(\Sigma(\Zset),\Omega(V), M_N(C^\infty (\Zset \times V)))$ the \vocab{matrix valued discrete-continuous bispectral context}.
We call a bispectral function 
\[
\Psi(k,x)\in M_N(C^\infty (\Zset \times V))
\]
for this context a \vocab{matrix valued discrete-continuous bispectral function}.
\ede
This bispectral context is \vocab{noncommutative} in the terminology of \deref{bisp context} if and only if the dimension of our matrices $N$ is greater than $1$.
We use capitalization to emphasize the fact that $\Psi(k,x)$ takes its values in the noncommutative algebra $M_N(\Rset)$ 
and write $\Psi(k,x)$ in place of $\Psi$ to emphasize the fact that it is a function.

The following is the simplest example of a discrete-continuous bispectral function:
\bex{ex1}
Let $V:=(0, + \infty)$ and $\Psi(k,x) := x^k I_N$. 
Recall that $I_N$ denotes the identity matrix of size $N \times N$.
Then $\Psi(k,x)$ is smooth on $\mathbb{Z}\times V$, has trivial left and right annihilators, and satisfies
$$\mathscr S_k \cdot\Psi(k,x) = \Psi(k,x)x\ \ \text{and}\ \ \Psi(k,x)\cdot \partial_xx = k\Psi(k,x),$$
so that $\Psi(k,x)$ is a discrete-continuous bispectral function.
The associated bispectral and Fourier algebras are
\begin{align}
&\mathcal B_k(\Psi) = M_N(\Rset)[\mathscr S_k ^{\pm 1}], 
&&\mathcal B_x(\Psi) = M_N(\Rset)[\partial_xx], 
\label{xy-bisp}
\\
&\mathcal F_k(\Psi) = M_N(\Rset)[k,\mathscr S_k^{\pm 1}],
&&\mathcal F_x(\Psi) = M_N(\Rset)[\partial_x,x^{\pm 1}].
\label{xy-Fourier}
\end{align}
The generalized Fourier map $b_\Psi : \mathcal F_k(\Psi) \to \mathcal F_x(\Psi)$ is the algebra isomorphism given by 
\begin{equation}
\label{xy-Fourier-map}
b_\Psi(k):=\partial_xx, \quad b_\Psi(\mathscr S_k) := x, \quad \mbox{and} \quad b_\Psi(A) :=A, \; \; \forall A \in M_N(\Rset).
\end{equation}
\eex

As a slightly more complicated example of a bispectral function, we can consider the classical discrete-continuous bispectral function associated with the Hermite polynomials.
\bex{ex2}
Let $V:=(0, + \infty)$ and
\[
\Psi(k,x) :=\sqrt{\frac{2^k}{\Gamma(k+1)}}\ _2F_0(-k/2;(1-k)/2;-1/x^2)x^ke^{-x^2/2}I_N.
\]
Then $\Psi(k,x)$ is smooth on $\mathbb Z\times V$, has trivial left and right annihilators, and satisfies
$$\frac{1}{\sqrt{2}}(\sqrt{k+1}\mathscr S_k + \sqrt{k}\mathscr S_k^{-1})\cdot\Psi(k,x) = \Psi(k,x)x\ \ \text{and}\ \ \Psi(k,x)\cdot (\partial_x^2 - x^2+1)  = -2 k\Psi(k,x),$$
so that $\Psi(k,x)$ is a discrete-continuous bispectral function.
The associated bispectral and Fourier algebras are
\begin{align*} 
&\mathcal B_k(\Psi) = M_N(\Rset)[\sqrt{k+1}\mathscr S_k + \sqrt{k}\mathscr S_k^{-1}],
&&\mathcal B_x(\Psi) = M_N(\Rset)[\partial_x^2+x^2],
\\
&\mathcal F_k(\Psi) = M_N(\Rset)[\sqrt{k+1}\mathscr S_k, \sqrt{k}\mathscr S_k^{-1}]
,
&&\mathcal F_x(\Psi) = [\partial_x,x].
\end{align*}
The generalized Fourier map $b_\Psi : \mathcal F_k(\Psi) \to \mathcal F_x(\Psi)$ is the algebra isomorphism given by 
$$
b_\Psi(\sqrt{k+1}\mathscr S_k) = \frac{1}{\sqrt{2}}(x-\partial_x), \; \; 
b_\Psi(\sqrt{k}\mathscr S_k^{-1}) = \frac{1}{\sqrt{2}}(x + \partial_x), \; \; 
b_\Psi(A) =A, \; \; \forall A \in M_N(\Rset).
$$
\eex

\subsection{Bifiltration}
For the discrete-continuous bispectral context in \deref{dc context}, the algebras $\Sigma(\Zset)$ and $\Omega(V)$ are filtered by the order and bandwidth of operators.
We can leverage this to obtain finite bi-filtrations of the left and right Fourier algebras.

For $D(x, \partial_x) \in \Omega(V)$, denote by $\ord(D(x, \partial_x))$ its order.

\bde{bifiltration}
Let $L (k, \mathscr S_k) = \sum_{j=-n}^n A_j(k)\mathscr S_k^j  \in \Sigma(\Zset)$ be a shift operator for some integer $n>0$.
The \vocab{bandwidth of $L (k, \mathscr S_k)$} is the value
$$\bw(L (k, \mathscr S_k)) := 2\max\{j\geq 0: A_{-j}(k)\neq 0\ \text{or}\ A_j(k)\neq 0\}.$$
If $\Psi$ is a discrete-continuous bispectral function and $L\in\mathcal F_k(\Psi)$ then the order of the differential operator $b_\Psi(L)$ is called the \vocab{co-order} of $L$, denoted $\cord(L)$.
Similarly we can define the \vocab{co-bandwidth} $\cbw( D)$ of $ D\in \mathcal F_x(\Psi)$ to be the bandwidth of $b_\Psi^{-1}( D)$.
\ede
The bandwidth and co-order of shift operators and differential operators
determine bi-filtrations of $\mathcal F_k(\Psi)$ and $\mathcal F_x(\Psi)$
\begin{align*}
\mathcal F_k^{\ell,m}(\Psi) &:= \{L\in \mathcal F_k(\Psi) : \bw(L) \leq \ell,\ \cord(L)\leq m\},\\
\mathcal F_x^{\ell,m}(\Psi) &:= \{ D\in \mathcal F_x(\Psi) : \ord( D) \leq \ell,\ \cbw( D)\leq m\}
\end{align*}
(which are $\Rset$-subspaces of $\mathcal F_k(\Psi)$ and $\mathcal F_x(\Psi)$, respectively)
and the generalized Fourier map restricts to an $\Rset$-linear isomorphism 
\begin{equation}
\label{lin-isom}
b_\Psi : \mathcal F_k^{\ell,m}(\Psi) \stackrel{\cong}{\longrightarrow} \mathcal F_x^{m,\ell}(\Psi).
\end{equation}

\bex{ex1 contd} Assume the setting of \exref{ex1}.
Applying Eqs. \eqref{xy-Fourier} and \eqref{xy-Fourier-map}, one obtains that the bifiltrations of the associated Fourier algebras is
\begin{align}
\mathcal F_k^{2\ell,2m}(\Psi) &= \Span\{k^j\mathscr S_k^iA : 0\leq j\leq 2m,\ |i| \leq \ell, \ A \in M_N(\mathbb R)\},
\label{xy-F1}
\\
\mathcal F_x^{2m,2\ell}(\Psi) &= \Span\{(\partial_xx)^jx^i A : 0\leq j\leq 2m,\ |i| \leq \ell, \ A \in M_N(\mathbb R)\},
\label{xy-F2}
\end{align}
and in particular, $\dim \mathcal F_k^{2\ell,2m}(\Psi) = (2\ell+1)(2m+1) N^2$. Furthermore, the linear isomorphism \eqref{lin-isom}
is given by
\begin{equation}
\label{xy-b-lm}
b_\Psi(k^j\mathscr S_k^i A ) = (\partial_xx)^j x^i A
\end{equation}
for all $0\leq j\leq 2m$, $|i| \leq \ell$, and $A \in M_N(\mathbb R)$.
\eex
\sectionnew{Orthogonal polynomials}
\label{3}
In this section we focus on {\em{classical}} discrete-continuous bispectral functions 
which are those that are eigenfunctions of a second order differential operator and a 
shift operator of bandwidth two. We describe a construction
for obtaining those from classical orthogonal polynomials on the real line.
\subsection{Basic definitions}
\bde{weight function}
Let $(x_0,x_1)\subseteq\Rset$ be an open interval which can be possibly infinite.
A \vocab{weight} of size $N$ on $(x_0,x_1)$ is a $M_N(\Rset)$-valued smooth function 
$W: (x_0,x_1)\rightarrow M_N(\Rset)$ with $W(x)$ a positive definite almost everywhere $N\times N$ matrix, with finite moments $\int_{x_0}^{x_1} |x|^n W(x) dx$ for all $n\geq 0$.
The interval $(x_0,x_1)$ is called the \vocab{support of $W$}.
\ede
A weight function defines an inner product on the space of $M_N(\Rset)$-valued polynomials defined by
$$\langle P(x),Q(x)\rangle_W := \int_{x_0}^{x_1} P(x)W(x)Q(x)^*|dx|.$$
By Gram-Schmidt, the above inner product defines a sequence of pairwise orthogonal polynomials.
\bde{orthogonal polynomials}
A \vocab{sequence of orthogonal polynomials} for a weight function $W(x)$ with support $(x_0,x_1)$ is a sequence of $M_N(\Rset)$-valued 
polynomials $\{P_n(x)\}_{n=0}^\infty$ which are pairwise-orthogonal with respect to $\langle\cdot,\cdot\rangle_W$ with $P_n(x)$ of degree $n$ for each $n$.
\ede
A sequence of polynomials is unique up to a choice of normalization.
The sequence is unique if we impose an additional condition, such as $P_n(x)$ is monic for all $n$.  The sequences that satisfy $\langle P_n(x), P_n(x)\rangle_W =I$ for all $n$ are called orthonormal.

It is easy to see that a sequence of orthogonal polynomials automatically satisfies a three-term recurrence relation.
\bpr{3-term recurrence} \cite{K}
Let $\{P_n(x)\}_{n=0}^\infty$ be a sequence of $M_N(\Rset)$-valued orthogonal polynomials for a weight function $W(x)$ supported on $(x_0,x_1)$.
Then there exist sequences $\{A_n\}$, $\{B_n\}$ and $\{C_n\}$ in $M_N(\Rset)$ satisfying
$$xP_n(x) = A_nP_n(x) + B_nP_n(x) + C_nP_{n-1}(x),\ \forall n\geq 1.$$
\epr

We are interested in special sequences of orthogonal polynomials which are eigenfunctions of a second-order differential operator
\bde{classical polynomials}
A sequence of $M_N(\Rset)$-valued orthogonal polynomials for a weight function $W(x)$ supported on $(x_0,x_1)$ are \vocab{eigenfunctions of a second-order differential equation} if for every $n$
$$P_n(x)''A_2(x) + P_n(x)'A_1(x) + P_n(x)A_0(x) = \Lambda_n P_n(x)$$
for some sequence of matrices $\{\Lambda_n\}$.
In the special case that $N=1$, the sequence is called a sequence of \vocab{classical orthogonal polynomials on the real line}.
\ede

In other words, classical orthogonal polynomials are polynomials that are simultaneously eigenfunctions of a second-order differential operator and a shift operator of bandwidth 2. 
The classical orthogonal polynomials for an interval $(x_0,x_1)\subseteq\Rset$ were classified by Bochner \cite{bochner1929sturm} and are up to affine transformation given by the Hermite, Laguerre, and Jacobi polynomials.
These examples of orthogonal polynomials, described in terms of hypergeometric functions, along with their weights, supports, and differential operators are listed in the table in Figure \ref{classical polynomials table}.
\begin{figure}[htp]
\label{classical polynomials table}
\begin{center}
{\setlength{\extrarowheight}{3pt}%
\begin{tabular}{|c|c|c|c|}\hline
Name & \text{weight} & \text{support} & \text{operator}\\\hline
Hermite   & $e^{-x^2}$          & $\Rset$      & $\partial_x^2-2\partial_xx$\\\hline
Laguerre  & $x^be^{-x}$         & $(0,\infty)$ & $\partial_x^2x+\partial_x(b+1-x)$\\\hline
Jacobi    & $(1-x)^a(1+x)^b$    & $(-1,1)$     & $\partial_x^2(1-x^2) + \partial_x(b-a-(b+a+2)x)$\\\hline
\end{tabular}
}
\\\vspace{0.2in}
{\setlength{\extrarowheight}{3pt}%
\begin{tabular}{|c|c|}\hline
Name & $p_n(x)$\\\hline
Hermite   & $_2F_0(-k/2; (1-k)/2; -1/x^2)(2x)^k$\\\hline
Laguerre  & $_1F_1(-k; a+1; x)$\\\hline
Jacobi    & $_2F_1(-k;1+a+b+k;a+1;(1-x)/2)$\\\hline
\end{tabular}
}
\end{center}
\caption{The classical orthogonal polynomials on $\Rset$.}
\end{figure}

Note that for Laguerre $b>-1$ and for Jacobi both $a>-1$ and $b>-1$. 

\subsection{Classical discrete-continuous bispectral functions}
We are interested in matrix valued discrete-continuous bispectral functions, so we define classical orthogonal polynomials in this context.
\bde{classical polys}
Let $V\subseteq\Rset$ be an open interval.
We define a \vocab{classical discrete-continuous bispectral function} to be a discrete-continuous bispectral function $\Psi: \Zset\times V\rightarrow M_N(\Rset)$ 
taking values in $\Rset I_N$ (i.e. the scalars in $M_N(\Rset)$) 
which has nonconstant elements $L\in \mathcal F_{k,sym}^{2,0}(\Psi)$ and $ D\in\mathcal F_{x,sym}^{0,2}(\Psi)$ with 
$\Rset I_N$-valued (scalar) coefficients, and nonsingular leading coefficients, and which satisfies the orthonormality condition
$$\Psi(k,x)\in L^2(V) \; \; \mbox{for} \; \; k\geq 0\ \ \text{and}\ \ \int_V \Psi(n,x)\Psi(m,x)^* dx = \delta_{m,n} \; \; \mbox{for} \; \;  m,n\geq 0.$$
\ede

Classical discrete-continuous bispectral functions are intimately linked with classical orthogonal polynomials.
In particular, if $\Psi(k,x)\Psi(0,x)^{-1}$ is a polynomial of degree $k$ for all integers $k\geq 0$, 
then these polynomials will be orthogonal with respect to the weight $W(x) = \Psi(0,x)\Psi(0,x)^*$ and will satisfy a second-order differential equation.
By Bochner's classification theorem \cite{bochner1929sturm}, the sequence $\Psi(k,x)\Psi(0,x)^{-1}$ will be an affine transformation of the Hermite, Laguerre, or Jacobi polynomials listed in the table above (up to a choice of normalization).
The associated shift operator then specializes to the three-term recurrence relation for the orthogonal polynomials.

Conversely, given any of the classical orthogonal polynomials from the table, 
the associated hypergeometric function $\psi(k,x)$ defines a classical discrete-continuous 
bispectral function $\Psi(k,x)$ on $V$ for $V$ the support of the associated weight $W(x)$ and
$$\Psi(k,x) = \theta(k)\psi(k,x)\rho(x)$$
for well-chosen normalizating functions $\theta(k)$ and $\rho(x)$ which symmetrizes the shift operator associated to the three-term recurrence relation of the given specialization of a hypergeometric function.

The values of these discrete-continuous bispectral functions, along with their associated shift and differential operators, are provided in the table in Figure \ref{classical discrete-continuous table}.
\begin{figure}
\label{classical discrete-continuous table}
\begin{center}
{\setlength{\extrarowheight}{3pt}%
\begin{tabular}{|c|c|}\hline
Name & $\Psi(k,x)$\\\hline
Hermite   & $\sqrt{\frac{2^k}{\Gamma(k+1)}}\ _2F_0(-k/2;(1-k)/2;-1/x^2)x^ke^{-x^2/2}I$\\\hline
Laguerre  & $\sqrt{\frac{\Gamma(k+a+1)}{\Gamma(k+1)}}\frac{1}{\Gamma(a+1)}\ _1F_1(-k; a+1; x)x^{a/2}e^{-x/2}$\\\hline
Jacobi    & $\sqrt{\frac{(2k+a+b+1)\Gamma(k+a+b+1)\Gamma(k+a+1)}{\Gamma(k+1)\Gamma(k+b+1)}}\ _2F_1(-k;1+a+b+k;a+1;(1-x)/2)(1-x)^{a/2}(1+x)^{b/2}$\\\hline
\end{tabular}
}
\\\vspace{0.2in}
{\setlength{\extrarowheight}{3pt}%
\begin{tabular}{|c|c|c|}\hline
Name & $ D$ & $L$\\\hline
Hermite & $\partial_x^2-x^2+1$ & $\sqrt{k/2}\mathscr S_k^{-1} + \sqrt{(k+1)/2}\mathscr S_k$\\\hline
Laguerre & $\partial_xx\partial_x + \frac{1}{2}-\frac{(a-x)^2}{4x}$ & \shortstack{$A(k)\mathscr S_k^{-1} + (2k+1+a)+A(k+1)\mathscr S_k$\\ $A(k) = -\sqrt{k(k+a)}$}\\\hline
Jacobi & \shortstack{$\partial_x(1-x^2)\partial_x + \frac{a^2}{2(x-1)}- \frac{b^2}{2(x+1)}$\\ $ + \frac{(a+b)(a+b+2)}{4}$} & \shortstack{$A(k)\mathscr S_k^{-1} + \frac{b^2-a^2}{(2k+2+a+b)(2k+a+b)}+A(k+1)\mathscr S_k$\\ $A(k) = \sqrt{\frac{(k+a+1)(k+b+1)k(k+a+b)}{(2k+a+b)^2-1}}\left(\frac{2}{2k+a+b}\right)$}\\\hline
\end{tabular}
}
\end{center}
\caption{Classical discrete-continuous bispectral functions.}
\end{figure}

\sectionnew{Behaviour of Fourier algebras under bispectral Darboux transformations}
\label{4}
In this section we prove sharp estimates on the sizes of the bifiltrarions of the Fourier algebras 
of all bispectral Darboux transformations from classical discrete-continuous bispectral functions.

\subsection{Estimates in the classical case}
For a matrix valued discrete-continuous bispectral function $\Psi(k,x)$ define the subspaces
\begin{align*}
\mathcal F_{k,sym}^{\ell,m}(\Psi) &:= \Span\{ L \in F_k^{\ell,m}(\Psi) : L^* = L, b_\Psi(L)^* = b_\Psi(L) \},
\\
\mathcal F_{x,sym}^{m,\ell}(\Psi) &:= \Span\{ D \in F_x^{m,\ell}(\Psi) : D^* = D, b_\Psi^{-1}(D)^* = b_\Psi(D) \}.
\end{align*}
It follows from \eqref{lin-isom} that the generalized Fourier map $b_\Psi$ maps $\mathcal F_{k,sym}^{\ell,m}(\Psi)$ bijectively onto 
$\mathcal F_{x,sym}^{m,\ell}(\Psi)$.
\bre{b-star}
Note that when $\Psi$ is a self-adjoint bispectral Darboux transformation of a classical discrete-continuous bispectral function, the generalized Fourier map 
$b_\Psi$ satisfies 
\[
b_\Psi(L^*) = b_\Psi(L)^* \quad \mbox{for all} \quad L \in \mathcal F_{k}^{\ell,m}(\Psi).
\] 
\ere

\ble{fundamental est lem}
Let $\Psi(k,x)$ be a classical matrix valued discrete-continuous bispectral function. Then
$$\dim \mathcal F_{k,sym}^{2\ell,2m}(\Psi)  = \dim \mathcal F_{x,sym}^{2m,2\ell}(\Psi) \geq (\ell + 1)(m+1) N^2.$$
\ele
\begin{proof}
Recall \eqref{lin-isom}. We will show the inequality for $ \dim \mathcal F_{x,sym}^{m,\ell}$.
By definition, there exist $ D\in\mathcal F_{x,sym}^{2,0}(\Psi)$ and $L\in\mathcal F_{k,sym}^{2,0}(\Psi)$ with
\begin{align*}
&\Psi(k,x)\cdot D = \Lambda(k)\Psi(k,x) \quad &&\mbox{for} \quad  D = (\partial_x^2d_2(x) + \partial_xd_1(x) + d_0(x))I \quad \mbox{and}
\\
&L\cdot\Psi(k,x) = \Psi(k,x) F(x) \quad &&\mbox{for} \quad L = (c(k)\mathscr S_k^{-1} + b(k) + a(k)\mathscr S_k)I.
\end{align*}
Since $\Psi(k,x)$ takes its values in $\Rset I_N$, the spaces $\mathcal F_{x,sym}^{2m,2\ell}(\Psi)$ are also $(M_N(\Rset),M_N(\Rset))$-bimodules.
Consequently we see that 
$$\mathcal F_{x,sym}^{2m,2\ell}(\Psi) \supseteq \Span \{A D^i F(x)^j +  F(x)^j  D^iA^*: 0\leq i\leq m,\ 0\leq j\leq \ell,\ A\in M_N(\Rset)\}.$$
By comparing orders and leading coefficients, we see that the right hand side is an $\Rset$-vector space of dimension at least $(\ell+1)(m+1)N^2$.
\end{proof}
The lower bound in \leref{fundamental est lem} is sharp as shown in next example.
\bex{cont-two ex1} Consider the setting of Examples \ref{eex1} and \ref{eex1 contd} with the classical discrete-continuous bispectral function $\Psi(k,x) = x^k I_N$. Using eqs. 
\eqref{xy-F1} and \eqref{xy-b-lm}, one easily shows that 
\[
\mathcal F_{x,sym}^{2m,2\ell}(\Psi) = \Span \{A (\partial_x x^2 \partial_x)^i x^j +  x^j (\partial_x x^2 \partial_x)^iA^*: 0\leq i\leq m,\ 0\leq j\leq \ell,\ A \in M_N(\Rset)\}.
\]
Therefore, in this case, $\dim F_{k,sym}^{2\ell,2m}(\Psi) =  (\ell + 1)(m+1)N^2.$
\eex
\subsection{Estimates for all bispectral Darboux transformations}
\bde{bispectral darboux degree}
Let $\wt\Psi(k,x)$ be a bispectral Darboux transformation of a discrete-continuous bispectral function $\Psi(k,x)\in M_N(C^\infty (\Zset \times V))$ with
\begin{align*}
\wt\Psi(k,x) &= F(k)^{-1}P(k, \mathscr S_k) \cdot \Psi(k,x) Q(x)^{-1}
\\
\Psi(k,x) &= \wt{P} (k, \mathscr S_k) \wt F(k)^{-1}\cdot \wt\Psi(k,x)\wt Q(x)^{-1}
\end{align*}
for some shift operators $P (k, \mathscr S_k),\wt{P}(k, \mathscr S_k) \in \mathcal F_k(\Psi)$ 
and functions $F(k),\wt F(k) : \Zset \to M_N(\Rset)^\times$, $Q(x),\wt Q(x) \in M_N(C^\infty (V))^\times$, cf. \deref{bispectral}.

We say that the bispectral Darboux transformation is \vocab{robust} if the right annihilator of $P$ in $\Sigma(\Zset)$ is trivial, 
the left annihilator of $\wt{P}$ in $\Sigma(\Zset)$ is trivial, the left annihilator of $b_\Psi(P)$ in $\Omega(V)$ is trivial, 
and the right annihilator of $b_\Psi(\wt{P})$ in $\Omega(V)$ is trivial
(recall the notation for the discrete-continuous bispectral context in \deref{dc context}).
\ede

Note that in the case $N=1$ every bispectral Darboux transformation is robust. 

Recall that the bispectral Darboux transformation is self-adjoint if 
$F(k)^*=\wt F(k)$, $Q(x)^* = \wt Q(x),$ $P^* = \wt{P}$ and $b_\Psi(P)^* = b_\Psi( \wt{P})$.
We call $(d_1,d_2)$ the \vocab{degree} of the self-adjoint bispectral Darboux transformation where 
$$d_1 = \bw(P) \quad \mbox{and} \quad d_2= \ord(b_\Psi(P)).$$

Our next result gives a crucial sharp estimate on the behaviour of the growth of Fourier algebras under Darboux 
transformations. This theorem is at the heart of linking matrix valued discrete-continuous bispectrality to the prolate spheroidal property.
\bth{fundamental est thm}
Let $\wt\Psi(k,x)$ be a robust, self-adjoint bispectral Darboux transformation of degree $(d_1,d_2)$ of a discrete-continuous bispectral function $\Psi(k,x)$
with values in $M_N(\Rset)$. Assume moreover that $\mathcal F_x^{*,0}(\Psi)$ contains a self-adjoint differential operator $C(x,\partial_x)$ of order $2d$ whose leading coefficient is scalar valued.  Then
\begin{align*}
\dim \mathcal F_x^{2\ell,2m}(\wt\Psi) &\geq \dim \mathcal F_x^{2\ell,2m-2d_2}(\Psi) + \dim \mathcal F_x^{2\ell-2d_1,2d_2-2d}(\Psi) + 1,
\\
\dim \mathcal F_{x,sym}^{2\ell,2m}(\wt\Psi) &\geq \dim \mathcal F_{x,sym}^{2\ell,2m-2d_2}(\Psi)  + \dim \mathcal F_{x,sym}^{2\ell-2d_{1},2d_2-2d}(\Psi) + 1.
\end{align*}
\eth
\begin{proof}
By assumption, 
$$\wt\Psi(k,x) = F(k)^{-1}P\cdot \Psi(k,x) Q(x)^{-1}\ \ \text{and} \ \ \Psi(k,x) = P^* (F(k)^{-1})^*\cdot \wt\Psi(k,x)(Q(x)^{-1})^*$$
for some $P\in \mathcal F_{k,sym}(\Psi)$, $F(k): \Zset\rightarrow M_N(\Rset)^\times$ and $Q(x)\in M_N(C^\infty(V))^\times$ with 
\[
d_1 = \bw(P) \quad \mbox{and} \quad d_2 = \ord(b_\Psi(P)).
\]
Let $T(x,\partial_x) := b_\Psi(P)$ and consider the vector space
\begin{align*}
\mathcal V_{\ell,m}
  & := F(k)\mathcal F_k(\Psi)F(k)^*\cap \mathcal F_k^{\ell,m}(\wt \Psi)\\
  & = \{F(k)L(x,\mathscr S_x) F(k)^*: L\in \mathcal F_k(\Psi),\ \bw(L)\leq \ell,\ \ord (Tb_\Psi(L)T^*) \leq m\}.
\end{align*}
Since our bispectral Darboux transformation is robust, $T$ has no nonzero right annihilator and $T^*$ has no nonzero left annihilator.  Therefore $\mathcal V_{\ell,m}$ is finite dimensional.

The left Fourier algebra of $\wt\Psi$ satisfies
$$\mathcal F_k^{2\ell,2m}(\wt\Psi)\supseteq \mathcal V_{2\ell,2m} + (F(k)^{-1}P) 
\mathcal F_k^{2\ell-2d_1,2d_2-2d}(\Psi)(F(k)^{-1}P)^* + \Rset I.$$
The sum on the right hand side may not be direct.  However, we have
$$\mathcal V_{\ell,m}\cap (F(k)^{-1}P) 
\mathcal F_k^{2\ell-2d_1,2d_2-2d}(\Psi)(F(k)^{-1}P)^* + \Rset I \subseteq \mathcal F_k^{2\ell,2d_2-2d}(\Psi).$$
Therefore,
\begin{align*}
\dim(\mathcal F_k^{\ell,m}(\wt\Psi))
  & \geq \dim(\mathcal V_{\ell,m}) + \dim\mathcal F_k^{2\ell-2d_1,2d_2-2d}(\Psi) + 1 - \dim(\mathcal F_k^{2\ell,2d_2-2d}(\Psi)).
\end{align*}

For a fixed $\ell$, the algebra $\mathcal F_x^{*,2\ell}(\Psi)$ is a finitely generated, torsion-free $\mathbb R[C]$-module.  
Hence it is a free module and we can choose a basis to write
$$\mathcal F_x^{*,2\ell}(\Psi) = \bigoplus_{i}\mathbb R[C(x,\partial_x)]A_i(x,\partial_x).$$
Since the leading coefficient of $C(x,\partial_x)$ is scalar valued, 
\[
\ord(TC^mA_iT^*) = \ord(TA_iT^*) + md.
\]
In fact,  we have for all $m > \max_i\ord(TA_iT^*)$,
$$\dim(\mathcal V_{2\ell,2m}) = \sum_{i} (1+\lfloor (2m-\ord(TA_iT^*))/2d\rfloor).$$
Likewise, for all $2m > \max_i \{ \ord(A_i) \}$,
$$\dim(\mathcal F_k^{2\ell,2m}(\Psi)) = \sum_{i}(1+\lfloor (2m-\ord(A_i))/2d\rfloor).$$
Let $\epsilon_i = (2d_2 + \ord(A_i) - \ord(TAT^*))/2d$ and note that $\epsilon_i\geq 0$ for all $i$. Then we have
\begin{align*}
\dim(\mathcal V_{2\ell,2m})
  & = \sum_{i} (1+\lfloor (2m-2d_2-\ord(A_i))/2d + \epsilon_i\rfloor)\\
  & \geq  \sum_{i} (1+\lfloor (2m-2d_2-\ord(A_i))/2d \rfloor)  + \sum_{i} \lfloor \epsilon_i\rfloor\\
  & \geq  \dim \mathcal F_k^{2\ell,2m-2d_2}(\Psi) + \dim \mathcal F_k^{2\ell,2d_2-2d}(\Psi).
\end{align*}

Therefore,
\begin{align*}
\dim(\mathcal F_k^{2\ell,2m}(\wt\Psi)) \geq  \dim(\mathcal F_k^{2\ell,2m-2d_2}) + \dim\mathcal F_k^{2\ell-2d_1,2d_2-2d}(\Psi) + 1.
\end{align*}
The same argument, replacing $\mathcal F_k$ with $\mathcal F_{k,sym}$ provides the second estimate.
\end{proof}

\bco{dimension after transformation}
Let $\wt\Psi(k,x)$ be a robust, self-adjoint bispectral Darboux transformation of degree $(d_1,d_2)$ of a classical discrete-continuous bispectral function $\Psi(x,y)$.
Then
$$\dim \mathcal F_{k,sym}^{2\ell,2m}(\wt\Psi) \geq 1+ \big( (\ell+1)(m+1) - d_1d_2 \big)N^2.$$
\eco
\begin{proof}
Since $\Psi(x,y)$ is a classical discrete-continuous bispectral function, there exists $C(x,\partial_x) \in \mathcal F_x^{2,0}(\Psi)$.
The estimate now follows from Lemma \ref{lfundamental est lem} and Theorem \ref{tfundamental est thm}.
\end{proof}
We note that the lower bound in the general case treated in \coref{dimension after transformation} differs by a constant from the lower 
bound in the classical case treated in \leref{fundamental est lem}. This will play a fundamental role in our approach to 
obtaining commutativity between integral and differential operators, and the discrete analog of this result.
\sectionnew{Adjoints and bilinear concomitants}
\label{5}
In this section we treat bilinear concomitants of differential and difference operators
which are used to construct self-adjoint operators in the Fourier algebras  
of matrix valued discrete-continuous bispectral functions.
For a general reference of adjointability and formal adjoints of differential operators, see \cite{coddington}.

In this section, as in \deref{dc context}, for a connected open interval $V \subseteq \Rset$, we will denote 
by $\Om(V)$ the opposite algebra of the algebra of $M_N(\Rset)$-valued smooth differential operators on $V$. 
\subsection{Adjoints of differential operators}
Differential operators have two possible notions of an adjoint, namely their adjoint as an unbounded linear operator on an inner product space and 
their \vocab{formal adjoint}, a natural anti-involution of the algebra of differential operators.
In general, these two adjoints may not be the same: the operator adjoint depends on the vector space of functions on which we act, which in turn depends on choices such as the domain.
When the linear operator acts as a differential operator on a dense subspace of its domain, the associated differential operator must be the formal adjoint.

\bde{continuous adj}
A differential operator $ D\in\Omega(V)$ is called \vocab{formally symmetric} if $ D^*= D$.
We say $ D\in\Omega(V)$ is \vocab{adjointable} with respect to an interval $(x_0,x_1)\subseteq V$ if the coefficients of $ D$ are smooth on $(x_0,x_1)$ and $\int_{x_0}^{x_1} (F(x)\cdot D)G(x)^* dx = \int_{x_0}^{x_1} F(x)(G(x)\cdot D)^*dx$ for all smooth, matrix valued functions $F(x),G(x)$ on $(x_0,x_1)$ with compact support.
If $ D$ is formally symmetric and adjointable with respect to $(x_0,x_1)$, then $ D$ is called \vocab{self-adjoint} with respect to $(x_0,x_1)$.
\ede

There is a nice algebraic description of the formally symmetric differential operators.
\ble{continuous formal adj lem}
Let $ D\in\Omega(V)$ be a differential operator of order $2m$ with $ D^* =  D$.
Then there exist $*$-symmetric matrix valued smooth functions $A_0(x),\dots, A_m(x)$ and $*$-skewsymmetric
matrix valued smooth functions $B_1(x),\dots, B_m(x)$ on $V$ satisfying
$$ D = \sum_{i=0}^m \partial_x^iA_i(x)\partial_x^i + \sum_{i=1}^m \{\partial_x^{2i-1},B_m(x)\},$$
where here $\{R,S\} = RS + SR$ denotes the anticommutator of operators $R$ and $S$.
\ele
The scalar case of this lemma can be traced back to \cite[Ch. 11]{coddington}.
\begin{proof}
This follows from a simple inductive argument on the order of $ D$.
\end{proof}

The difference between the formal adjoint and the operator adjoint is captured by its bilinear concomitant.

\bde{continuous conco} Consider an $M_N(\Rset)$-valued differential operator with smooth coefficients:
$$ D = \sum_{j=0}^m \partial_x^j A_j(x) \in \Omega(V).$$
The {\em{bilinear concomitant}} of $ D$ is the bilinear form $\mathcal{C}_{ D}(\cdot,\cdot; p)$ defined on smooth functions $F,G \in M_N(C^\infty(V))$ by
\begin{align*}
\mathcal{C}_D(F,G;x)
  & := \sum_{j=1}^m \sum_{i=0}^{j-1} (-1)^i F^{(j-1-i)}(x)(G(x) A_j(x)^*)^{(i)*}\\
  & = \sum_{j=1}^m \sum_{i=0}^{j-1}\sum_{\ell=0}^i\binom{i}{\ell} (-1)^i F^{(j-1-i)}(x)A_j^{(i-\ell)}(x)G(x)^{(\ell)*}.
\end{align*}
\ede

Equivalently, for the smooth matrix valued function $C_{ D}(x) \in M_m(C^\infty(V))$ 
whose $n,\ell$-th entry is given by
\begin{equation}
\label{C-matr}
\mathcal C_{ D}(x)_{n,\ell} := \sum_{j=\ell+n+1}^{m}\binom{j-n}{\ell-1} (-1)^{j-n} A_j(x)^{(j+1-n-\ell)},
\end{equation}
the bilinear concomitant may be expressed as
$$\mathcal{C}_{ D}(F,G;x) = [F(x)\ F'(x)\ \dots\ F^{(m)}(x)] C_{ D}(x)[G(x)\ G'(x)\ \dots\ G^{(m)}(x)]^{*}.$$
The concomitant $\mathcal C_{ D}$ is $M_N(\Rset)$-sesquilinear in the sense that it is $\Rset$-bilinear and
$$
\mathcal C_{ D}(CF,G; z) = C\mathcal C_{ D}(F,G; z)\ \ \text{and}\ \ 
\mathcal C_{ D}(F,CG; z) = \mathcal C_{ D}(F,G; z)C^*,\ \ \forall\ C\in M_N(\Rset).
$$

As the next proposition shows, the operator adjoint is linked with the formal adjoint via this concomitant.
Its proof is standard and is left to the reader. The scalar case of the proposition 
appears in \cite[Ch. 11]{coddington}.
\bpr{integral biconcomitant}\label{integral biconcomitant}
Suppose $(x_0,x_1)\subseteq V$ is a possibly infinite open interval.
Let $ D\in \Omega(V)$ and let $F(x),G(x)$ be two smooth, matrix valued functions on $(x_0,x_1)$ with compact support.  Then
$$\int_{x_0}^{x_1}\left[(F(x)\cdot D)G(x)^*-F(x)(G(x)\cdot  D^*)^*\right] dx = \mathcal{C}_{ D}(F,G;x_1) - \mathcal{C}_{ D}(F,G;x_0),$$
with the convention that $\mathcal{C}_D(F,G;\pm\infty) = 0$.
\epr

\subsection{Adjoints of shift operators}
The connection between adjoints of differential operators and concomitants is paralleled in the setting of difference operators.
\bde{discrete adj}
A shift operator $L (k, \mathscr S_k) \in\Sigma(\Zset)$ is called \vocab{formally symmetric} if 
\[
L(k, \mathscr S_k)^*=L(k, \mathscr S_k).
\]
We say $L\in\Sigma(\Zset)$ is \vocab{adjointable} with respect to a set $I\subseteq \Zset$ if
$$\sum_{k\in I}  (L (k, \mathscr S_k) \cdot F(k))^*G(k)^* = \sum_{k\in I} F(k)^*(L (k, \mathscr S_k) \cdot G(k))^*$$
for all $M_N(\Rset)$-valued functions $F(k),G(k)$ on $\Zset$, vanishing at all but finitely many points of $I$.
If $L(k, \mathscr S_k)$ is formally symmetric and adjointable with respect to $I$, then $L(k, \mathscr S_k)$
 is called \vocab{self-adjoint} with respect to $I$.
\ede
In other words, a shift operator is self-adjoint if and only if it is self-adjoint as a linear operator.
Conversely, being formally symmetric is an algebraic condition which is characterized by the following lemma.
\ble{discrete formal adj lem}
Let $L(k, \mathscr S_k) \in \Sigma(\Zset)$ be a shift operator with bandwidth $\ell$ and $L(k, \mathscr S_k) ^* = L(k, \mathscr S_k) $.
Then there exist $*$-symmetric matrix valued functions $A_0(k),\dots, A_\ell(k)$ on $\Zset$ and $*$-skewsymmetric 
matrix valued functions $B_1(k),\dots, B_\ell(k)$ on $\Zset$ satisfying
$$L(k, \mathscr S_k)  = A_0(k) + \sum_{i=1}^\ell (A_i(k-i)\mathscr S_k^{-i} + A_i(k)\mathscr S_k^i) 
+ \sum_{j=1}^\ell (B_j(k-j)\mathscr S_k^{-j} - B_j(k)\mathscr S_k^j).$$
All shift operators of this form are formally symmetric.
\ele
\begin{proof}
This follows from a simple inductive argument on the bandwidth of $L(k, \mathscr S_k)$.
\end{proof}

The connection between the operator adjoint and the formal adjoint is captured by a discrete analog of the continuous bilinear concomitant defined in the previous 
subsection.
For any integer $\ell>0$ and functions $F,G: \Zset\rightarrow M_N(\Rset)$ define an $M_N(\Rset)$-valued smooth function on $\Zset$ by
$$\mathcal C_{\mathscr S^\ell}(F,G; z) := \sum_{i=1}^{\ell} F(z+i)^*G(z+i-\ell).$$
More generally, for an arbitrary shift operator $L(k, \mathscr S_k) = \sum_{n=-\ell}^\ell A_n(k)\mathscr S_k^n\in \Sigma(\Zset)$ define
\begin{equation}\label{discrete concomitant}
B_{L}(F,G; z) := \sum_{n=1}^\ell \left(\mathcal C_{\mathscr S^n}(F,A_n^*G; z) - \mathcal C_{\mathscr S^n}(A_{-n}^*G,F;z)^*\right).
\end{equation}
Then $B_{L}$ defines a map which is $\Rset$-bilinear and $M_N(\Rset)$-sesquilinear in the sense that
$$
\mathcal C_{L}(FA,G; z) = A^*\mathcal C_{L}(F,G; z)\ \ \text{and}\ \ 
\mathcal C_{L}(F,GA; z) = \mathcal C_{L}(F,G; z)A,\ \ \forall\ A\in M_N(\Rset).
$$
\bde{discrete conco}
For each $L (k, \mathscr S_k) \in \Sigma(\Zset)$, the $M_N(\Rset)$-sesquilinear map $\mathcal C_{L}$ 
defined above is called the \vocab{bilinear concomitant of $L (k, \mathscr S_k)$}.
\ede

The bilinear concomitant of shift operators captures the boundary behavior of adjunction of differential operators in a way similar 
to Proposition \ref{integral biconcomitant} as stated in the next proposition. Its proof is simple and is left to the reader. 
\bpr{discrete integral adjoint}
Let $L (k, \mathscr S_k) \in \Sigma(\Zset)$ and let $F,G: \Zset\rightarrow M_N(\Rset)$ be zero at all but finitely many values of $\Zset$.
Then
$$\sum_{k=m}^n \left[(L\cdot F)(k)^*G(k) - F(k)^*(L^*\cdot G)(k)\right] = \mathcal C_{L}(F,G;n)-\mathcal C_{L}(F,G;m-1),$$
with the convention that $\mathcal C_{L}(F,G;\pm\infty) = 0$.
\epr

\subsection{Bi-self-adjoint operators}
The generalized Fourier map associated with a discrete-continuous bispectral function $\Psi(k,x)$ directly connects an algebra of difference operators with an algebra of differential operators.
In the next section, we will show that differential operators commuting with integral operators will emerge from self-adjoint difference operators whose images under the generalized Fourier map are self-adjoint differential operators.

\bde{bisymmetric}
Fix an $M_N(\Rset)$-valued discrete-continuous bispectral function $\Psi(k,x)$, and let $(x_0,x_1)\subseteq V$ be a subinterval and $I\subseteq \Zset$.
\begin{enumerate}
\item We call a shift operator $L\in \mathcal F_k(\Psi)$ 
\vocab{formally bisymmetric} if $L^*=L$ and $b_\Psi(L)^* = b_\Psi(L)$.
\item We say that $L\in \mathcal F_k(\Psi)$ is \vocab{bi-self-adjoint with respect to the pair $(I,(x_0,x_1))$} 
if $L$ is self-adjoint with respect to $I$ and $b_\Psi(L)$ is self-adjoint with respect to $(x_0,x_1)$.
\end{enumerate}
\ede

When $\Psi(k,x)$ is self-adjoint bispectral Darboux transformation of a classical discrete-continuous bispectral function, the generalized Fourier map preserves adjoints.
In this situation, a difference operator being formally bisymmetric is no different than being formally symmetric, or from the associated differential operator being formally symmetric.
However, the existence of bi-self-adjoint operators is more subtle, since it requires the simultaneous vanishing of both discrete and continuous bilinear concomitants.
The next theorem provides a sufficient condition for the existence of these operators.

\bth{bisymmetric existence}
Suppose that $\Psi \in M_N( C^\infty( \Zset\times V))$ is a discrete-continuous bispectral function and that
$$\dim \mathcal F_{x,sym}^{\ell,m}(\Psi) \geq (\ell+1)(m+1)N^2  - c$$
for some $c\in\Zset$.
Let $I = \{k_0,k_0+1,\dots,k_1\}\subseteq \Zset$ and $(x_0,x_1)\subseteq V$, and assume that the concomitant of every operator in $\mathcal F_k(\Psi)$ vanishes at $k_0$, and the concomitant of every operator in $\mathcal F_x(\Psi)$ vanishes at $x_0$.
Then there exists a nonconstant $L (k, \mathscr S_k) \in \mathcal F_{k,sym}(\Psi)$ such that $L (k, \mathscr S_k)$ 
is bi-self-adjoint with respect to $(I,(x_0,x_1))$.
Furthermore, $L (k, \mathscr S_k)$ may be taken to be in $\mathcal F_{k,sym}^{2\ell,2\ell}(\Psi)$ for $\ell N^2 > c$.
\eth
\begin{proof}
The subspace of $\mathcal F_{k,sym}^{2\ell,2\ell}(\Psi)$ consisting of difference operators which are self-adjoint with respect to $I$ is
$$\mathcal U = \left\lbrace L\in \mathcal F_k(\Psi): \mathcal C_{L}(F,G;k_1) = 0\right\rbrace.$$
Using \leref{discrete formal adj lem}, a formally symmetric difference operator 
$$L = A_0(k) + \sum_{i=1}^\ell (A_i(k-i)\mathscr S_k^{-i} + A_i(k)\mathscr S_k^i) + \sum_{j=1}^\ell (B_j(k-j)\mathscr S_k^{-j} - B_j(k)\mathscr S_k^j)$$
has concomitant vanishing at $k_1$ if and only if $A_i(k_1+j) = 0$ and $B_i(k_1+j)=0$ for all $1\leq j\leq i$.
Using the fact that $A_i(k) \in M_N(\Rset)$ is symmetric and $B_i(k_1) \in M_N(\Rset)$ is skewsymmetric, 
we see that $\mathcal U$ is a subspace of codimension at most $N^2\ell(\ell+1)/2$.

The subspace of $\mathcal F_{x,sym}^{2\ell,2\ell}(\Psi)$ consisting of difference operators which are self-adjoint with respect to $(x_0,x_1)$ is
$$\mathcal V = \left\lbrace D\in \mathcal F_x(\Psi): \mathcal C_{D}(F,G;x_1) = 0\right\rbrace.$$
Using \leref{continuous formal adj lem}, a formally symmetric differential operator
$$ D = \sum_{i=0}^m \partial_x^iA_i(x)\partial_x^i + \sum_{i=1}^m \{\partial_x^{2i-1},B_m(x)\},$$
has concomitant vanishing at $x_1$ if and only if $A_i^{(j)}(x_1) = 0$ and $A_i^{(j)}(x_1) = 0$ for all $0\leq j < i$, and likewise for $B_i^{(j)}(x)$.
Using the fact that $A_i(x)$ is symmetric and $B_i(x)$ is skew-symmetric, we see that $\mathcal V$ is also a subspace of codimension at most $N^2\ell(\ell+1)/2$.

Choose $\ell$ such that $N^2\ell>c$.
The bi-self-adjoint operators lie in the intersection $\mathcal U\cap \mathcal V$, 
which has codimension at most $N^2\ell(\ell+1)$.  Since $\dim \mathcal F_{x,sym}^{2\ell,2\ell}(\Psi) \geq N^2(\ell+1)^2 - c$, we see that $\mathcal U\cap\mathcal V$ has dimension at least $N^2(\ell + 1) -c$, which is greater than $N^2$.
Consequently it must contain a nonconstant bi-self-adjoint operator.
\end{proof}

\sectionnew{Main theorem}
\label{6}
In this section, we prove a generalization of the Main Theorem stated in the introduction.
Throughout this section, ${\wt\Psi}$ will be a discrete-continuous bispectral function ${\wt\Psi}: \Zset\times V\rightarrow M_N(\Rset)$.
Using it, we define the integral operator
\begin{equation}\label{eqn:cont integral repeated}
T_{\wt\Psi}: F(y)\mapsto \int_{x_0}^{x_1} F(x)K(x,y) dx,\ \ K(x,y) := \sum_{k=k_0}^{k_1} {\wt\Psi}(k, x)^*{\wt\Psi}(k,y),
\end{equation}
and the discrete integral operator
\begin{equation}\label{eqn:disc integral repeated}
S_{\wt\Psi}: F(m)\mapsto \sum_{k=k_0}^{k_1} J(m,k)F(k),\ \ J(m,n) := \int_{x_0}^{x_1} {\wt\Psi}(m, y) {\wt\Psi}(k,y)^* dy,
\end{equation}
where here $(x_0,x_1)\subseteq V$, and $*$ is the matrix transposition.

\subsection{General statement}
The connection between bi-self-adjoint differential or difference operators and operators commuting with $T_{\wt\Psi}$ or $S_{\wt\Psi}$ is established in the following theorem.
\bth{basic}
Suppose $L(k,\mathscr S_k)\in \mathcal F_k({\wt\Psi})$ is bi-self-adjoint with respect to $(I,\Gamma)$, 
and let $R(x,\partial_x) = b_{\wt\Psi}(L)$.  Then for any smooth, compactly supported function $F: V\rightarrow M_N(\Rset)$,
\begin{equation}
\label{commutation1}
T_{\wt\Psi}( F\cdot R(x,\partial_x))(x)= (T_{\wt\Psi}(F)\cdot R)(x).
\end{equation}
Likewise, for every function $G: \Zset\rightarrow M_N(\Rset)$,
\begin{equation}
\label{commutation2}
S_{\wt\Psi}( L(k,\partial_k)\cdot G)(k)= (L(k,\partial_k)\cdot S_{\wt\Psi}(G)\cdot R)(k).
\end{equation}
\eth
\begin{proof}
Let $F$ be a smooth function on $\Gamma$ with compact support. Then we have
\begin{align*}
T_{\wt\Psi}( F\cdot R(x,\partial_x))(y)
  & = \int_\Gamma \left(F\cdot R\right)(x)\;K(x,y) dx\\
  & = \int_\Gamma \left(F\cdot R\right)(x)\;\sum_{k=k_0}^{k_1} {\wt\Psi}(k, x)^*{\wt\Psi}(k,y) dx\\
  & \stackrel{\text{I}}{=} \int_\Gamma F(x)\sum_{k=k_0}^{k_1} \left({\wt\Psi}(k, x)\cdot R\right)^*{\wt\Psi}(k,y) dx\\
  & \stackrel{\text{II}}{=} \int_\Gamma F(x)\sum_{k=k_0}^{k_1} \left(L\cdot{\wt\Psi}(k, x)\right)^*{\wt\Psi}(k,y) dx\\
  & \stackrel{\text{I}}{=} \int_\Gamma F(x)\sum_{k=k_0}^{k_1} {\wt\Psi}(k, x)^*\left(L\cdot{\wt\Psi}(k,y) \right)dx\\
  & \stackrel{\text{II}}{=} \int_\Gamma F(x)\sum_{k=k_0}^{k_1} {\wt\Psi}(k, x)^*\left({\wt\Psi}(k,y)\cdot R \right)dx\\
  & = (T_{\wt\Psi}(F)\cdot R)(y),
\end{align*}
where in this chain of equalities, we used in (I) that $L$ is bi-self-adjoint and in (II) that $L\in \mathcal F_x({\wt\Psi})$.
This proves \eqref{commutation1}.

A similar argument proves that $L$ commutes with $S_{\wt\Psi}$ in the sense of Eq. \eqref{commutation2}.
\end{proof}

As a direct consequence of Corollary \ref{cdimension after transformation}, and Theorems \ref{tbisymmetric existence} and \ref{tbasic}, 
we have the following result.
\bth{gral}
Let ${\wt\Psi}(x,y)$ be a robust, self-adjoint bispectral Darboux transformation of a classical discrete-continuous bispectral function 
$\Psi(x,y)$ of degree $(d_1,d_2)$ supported on $V$. Let $k_0 = 0$ and $x_0$ be the left endpoint of the interval $V$.
Then the following hold:
\begin{enumerate}
\item There exists a differential operator  $R(x,\partial_x) \in F_{x,sym}^{2d_1d_2,2d_1d_2}({\wt\Psi})$ 
commuting with the continuous integral operator $T_{{\wt\Psi}}$.
\item The shift operator $L (k,\mathscr S_k) : = b_{ {\wt\Psi}}^{-1} (R) \in \mathcal F_{k,sym}^{2d_1d_2,2d_1d_2}({\wt\Psi})$ 
commutes with the discrete integral operator $S_{\wt\Psi}$;
\end{enumerate}
\eth

The differential operator in part (1) of the theorem is obtained by imposing concomitant constraints on the operators of 
$F_{x,sym}^{2d_1d_2,2d_1d_2}({\wt\Psi})$ and solving the linear system of equations on its coefficients. 
Theorems \ref{tfundamental est thm} and \ref{tbisymmetric existence} guarantee that this system will have a nontrivial solution.
The shift operator in part (2) of the theorem is obtained by applying the generalized Fourier map to the first operator.


\sectionnew{Examples}\label{7}
This section contains examples illustrating the different features of the \thref{gral} and 
the power of its applications. We focus on part (1) of the theorem 
constructing a commuting differential operator for the integral operator $T_{\Psi}$. 
In all cases, one can construct a commuting shift operator for the corresponding discrete
integral operator $T_{\wt{\Psi}}$ by applying the (inverse) of the generalized Fourier
map for $\Psi$. We leave the details of part (2) to the reader. 
 
\subsection{Classical examples}

To begin, we will derive differential operators commuting with integral operators whose kernels are defined in terms of 
Christoffel--Darboux kernels of classical orthogonal polynomials.
Up to conjugation, these are precisely the differential operators obtained in \cite{grunbaum1}.
However, even in this special case, the commuting operators are obtained in an 
intrinsic fashion using the generalzied Fourier algebras of discrete-continuous bispectral functions.

Let $\Psi(k,x)$ be a classical discrete-continuous bispectral function with
$$L(k,\mathscr S_k)\cdot\Psi(k,x) = \Psi(k,x)x\quad\text{and}\quad \Psi(k,x)\cdot D(x,\partial_x) = \lambda(k)\Psi(k,x)$$
for some polynomial $\lambda(k)$ and operators
$$L(k,\mathscr S_k) =  A(k)\mathscr S_k^{-1} + B(k) + A(k+1)\mathscr S_k$$
and
$$D(x,\partial_x) = \partial_xp(x)\partial_x + q(x).$$
Thus the generalized bispectral map satisfies $b_\Psi(L) = x$ and $b_\Psi(\lambda(k)) = D$.

Now let $(x_0,x_1)$ be the (possibly infinite) subinterval of $\mathbb R$ of the support of the sequence of classical orthogonal polynomials corresponding to $\Psi$.
Based on the theory developed in the previous sections, we define an integral operator $T_\Psi$ by
$$T_\Psi: f(z)\mapsto \int_{x_0}^t f(x)K(x,z) dx,\ \ K(x,z) := \sum_{k=0}^n \Psi(k, x)\Psi(k,z).$$

To find a differential operator commuting with $T_\Psi$, our theory prompts us to try to find the symmetric elements of the left and right Fourier algebra of $\Psi$ of low order.
Already, the calculations above show that $D(x,\partial_x)$ and $x$ are both in $\mathcal F_{x,sym}(\Psi)$, while $L(k,\mathscr S_k)$ and $\lambda(k)$ are in $\mathcal F_{k,sym}(\Psi)$.

The anticommutator of a pair of symmetric operators is also symmetric, so we can use it to obtain more symmetric elements.
In particular,
$$\{D,x\} = 2\partial_xxp(x)\partial_x + p'(x) + 2xq(x)$$
is another element of $\mathcal F_{x,sym}(\Psi)$ and 
$$\{L,\lambda(k)\} = A(k)(\lambda(k)+\lambda(k-1))\mathscr S_k^{-1} + 2\lambda(k)B(k) + A(k+1)(\lambda(k)+\lambda(k+1))\mathscr S_k$$
is another element of $\mathcal F_{k,sym}(\Psi)$.
Note, since $b_\Psi$ is an algebra isomorphism, it must preserve the anticommutator, i.e. $b_\Psi(\{L,\lambda(k)\}) = \{D,y\}$.
As a consequence, the order two component of the symmetric portion of the Fourier algebra $\mathcal F_{k,sym}^{2,2}(\Psi)$ has dimension at least $8$.
By \thref{bisymmetric existence} combined with \thref{basic}, it follows that $\mathcal F_{x,sym}^{2,2}(\Psi)$ will contain a non-constant differential operator commuting with $T_\Psi$.

To find the commuting operator, we start with a generic linear combination of the non-constant elements in $\mathcal F_{x,sym}^{2,2}(\Psi)$ we have described so far
\begin{align*}
L(k,\mathscr S_k)
  & = A(k)(\lambda(k)+\lambda(k-1)+\beta)\mathscr S_k^{-1} + (\beta+2\lambda(k)) B(k)+\alpha\lambda(k)\\
  & + A(k+1)(\lambda(k)+\lambda(k+1)+\beta)\mathscr S_k
\end{align*}
Bx linearity, its image $R(x,\partial_x) = b_\Psi(L)$ is
$$R(x,\partial_x) = \partial_x(2x+\alpha)p(x)\partial_x + p'(x) + (\alpha+2x) q(x) +\beta x.$$

\begin{figure}[ht]
\begin{center}
{\setlength{\extrarowheight}{3pt}%
\begin{tabular}{|c|c|}\hline
Name & $R(x,\partial_x)$ \\\hline
Hermite & $\partial_x2(x-t)\partial_x + 2(x-t)(1-x^2) + 2(2n+1)x$\\\hline
Laguerre & $\partial_x2x(x-t)\partial_x + 1 - (x-t)(1-(a-x)^2/2x) + (2n+1)x$\\\hline
Jacobi & $\partial_x2(1-x^2)(x-t)\partial_x - 2x + 2(x-t)\left(\frac{a^2}{2(x-1)}- \frac{b^2}{2(x+1)}+ \frac{(a+b)(a+b+2)}{4}\right)$\\
 & $+ (2(n+1)^2 + (a+b)(2n+1))x$\\\hline
\end{tabular}
}
\end{center}
\caption{Differential operators commuting with integral operators for the classical discrete-continuous bispectral functions.}
\label{classical continuous commuting}
\end{figure}

The operator that we are searching for is exactly the one where the discrete concomitant \eqref{discrete concomitant} of $L(k,\mathscr S_k)$ vanishes at $n$ and where the continuous concomitant of $D(x,\partial_x)$ vanishes at $t$.
The vanishing condition of the concomitants results in the equations $\lambda(n)+\lambda(n+1) + \beta = 0$ and $2t+\alpha = 0$.
Therefore the differential operator 
$$R(x,\partial_x) = \partial_x2(x-t)p(x)\partial_x + p'(x) + 2(x-t)q(x) - (\lambda(n)+\lambda(n+1))x$$
commutes with the continuous integral operator $T_\Psi$.
The table in Figure \ref{classical continuous commuting} gives explicit expressions for the commuting operators in each of the basic cases.

The integral operator we have defined only makes sense for integer values of $n>0$.
Similar integral operators can be define for $n<0$, for which one can obtain commuting differential operators.

\subsection{A Darboux transformation example}

Self-adjoint bispectral Darboux transformations of classical discrete-continuous bispectral functions provide more complicated examples of integral operators commuting with differential operators, at the cost of the associated operators being of higher order.
To demonstrate this, we consider a specific bispectral Darboux transformation of the Laguerre-type discrete-continuous bispectral function.

Let $\Psi(k,x)$ be the classical discrete-continuous bispectral function of Laguerre type defined above and
$$D(x,\partial_x) = \partial_xx\partial_x - \frac{(a-x)^2}{4x} + \frac{1}{2},\quad\text{and}$$
$$L(k,\mathscr S_k) = -\sqrt{k(k+a)}\mathscr S_k^{-1} + (2k + a + 1) - \sqrt{(k+1)(k+a+1)}\mathscr S_k,$$
be the associated formally symmetric differential and shift operators.

Bispectral Darboux transformations come from rational factorizations of polynomials in the operator $D(x,\partial_x)$.
For one explicit example, the differential operator
$$Q(x,\partial_x) := \partial_x^2xq(x) + \partial_xq(0) - \left(\frac{x}{4} - \frac{2\lambda+a}{2} + \frac{a^2}{4x}\right)q(x),$$
for $q(x) = 2\lambda +2a+\sqrt{\frac{\lambda+a}{\lambda}}(2\lambda + a-x)$ defines a symmetric factorization of a quadratic polynomial in the Laguerre operator, via
$$Q(x,\partial_x)\frac{1}{q(x)^2}Q(x,\partial_x)^* = (D(x,\partial_x)+\lambda)(D(x,\partial_x) + \lambda-1).$$
Moreover, $Q(x,\partial_x)$ is an element of the right Fourier algebra of $\Psi(k,x)$, i.e. there exists a shift operator $P(k,\mathscr S_k)$ satisfying $P(k,\mathscr S_k)\cdot\Psi(k,x) = \Psi(k,x)\cdot Q(x,\partial_x)$.  
The explicit value of $P(k,\mathscr S_k)$ is obtained from leveraging the fact 
that the generalized Fourier map is an algebra isomorphism. 
In particular, we find
\begin{align*}
P(k,\mathscr S_k)
  & = b_{\Psi}^{-1}(Q(x,\partial_x))\\
  & = -\beta(k-\lambda)\sqrt{(k+1)(k+a+1)}\mathscr S_k - \beta(k+1-\lambda)\sqrt{k(k+a)}\mathscr S_k^{-1}\\
  & + 2\beta k(k+1) - 2\lambda\beta(\beta+2)k - \lambda \beta(\beta+2) + 2\lambda\beta a/(\beta-1),
\end{align*}
for $\beta = \sqrt{(\lambda+a)/\lambda}$.
The shift operators $P(k,\mathscr S_k)$ define the factorizations
$$P(k,\mathscr S_k)^*\frac{1}{p(k)^2} P(k,\mathscr S_k) = q(L(k,\mathscr S_k))^2.$$
for $p(k) = \sqrt{(\lambda-k)(\lambda-k-1)}$.  
Thus the functions
$$\wt{\Psi}(k,x) := \Psi(k,x)\cdot Q(x,\partial_x)\frac{1}{p(k)q(x)}$$
are self-adjoint bispectral Darboux transformations of $\Psi(k,x)$.

A self-adjoint bispectral Darboux transformation always leads to another self-adjoint bispectral function (see \deref{bisp context}).
In our situation, $\wt\Psi(k,x)$ are discrete-continuous bispectral functions and satisfy the differential equation
$$\wt{\Psi}(k,x)\cdot \frac{1}{q(x)}Q(x,\partial_x)^*Q(x,\partial_x)\frac{1}{q(x)} = p(k)^2\wt{\Psi}(k,x)$$
along with the difference equation
$$\frac{1}{p(k)} P(k,\mathscr S_k) P(k,\mathscr S_k)^*\frac{1}{p(k)}\cdot\wt{\Psi}(k,x) = \wt{\Psi}(k,x)q(x)^2.$$
Note that both the shift and differential operators above are formally symmetric.

Just as we have done previously, we use $\wt\Psi(k,x)$ to define an integral operator
$$T_{\wt\Psi}: f(y)\mapsto \int_{0}^t f(x)K(x,y) dx,\ \ K(x,y) := \sum_{k=0}^n \wt\Psi(k, x)\wt\Psi(k,y).$$
We will construct a differential operator which commutes with this integral operator.
In order to do so, we must construct more examples of formally symmetric operators in the left and right Fourier algebras.

Our main tools for constructing more formally symmetric operators in the Fourier algebra is the pair of reciprocal relations
\begin{align*}
\wt{\Psi}(k,x)\cdot \frac{1}{q(x)}Q(x,\partial_x)^*B(x,\partial_x)Q(x,\partial_x)\frac{1}{q(x)} = p(k) A(k,\mathscr S_k)p(k)\cdot\wt\Psi(k,x),\\
\frac{1}{p(k)} P(k,\mathscr S_k)^* 
A(k,\mathscr S_k) P(k,\mathscr S_k)\frac{1}{p(k)}\cdot\wt\Psi(k,x) = \wt\Psi(k,x)\cdot q(x)B(x,\partial_x)q(x),
\end{align*}
for all $A(k,\mathscr S_k)\in \mathcal F_k(\Psi)$ and $B(x,\partial_x)\in\mathcal F_x(\Psi)$ with $b_\Psi(A) = B$.
These allow us to build operators in $\mathcal F_k(\wt\Psi)$ and $\mathcal F_x(\wt\Psi)$ from operators in $\mathcal F_k(\Psi)$ and $\mathcal F_x(\Psi)$.
In particular, formal symmetricity is preserved by the construction.

To obtain our commuting differential operator, we use the previous paragraph to construct formally symmetric differential operators in 
$\mathcal F_x(\wt\Psi)$ of order $\leq 4$, whose preimages in $\mathcal F_k(\wt\Psi)$ under the generalized Fourier map have bandwidth $\leq 4$.
By showing that $\mathcal F_{x,sym}^{4,4}(\wt\Psi)$ has sufficiently high dimension, we can guarantee that it contains a differential operator commuting with our integral operator.

For starters, we can take $A(k,\mathscr S_k) = L(k,\mathscr S_k)^j$ and $B(x,\partial_x) = x^j$.
This gives us three linearly independent operators 
$$R_{j+1}(x,\partial_x) = \frac{1}{q(x)}Q(x,\partial_x)^*x^jQ(x,\partial_x)\frac{1}{q(x)},\quad\text{for}\ j=0,1,2.$$
Analogously, we can take $A(k,\mathscr S_k) = (-k)^j$ and $B(x,\partial_x) = D(x,\partial_x)^j$.
This gives us two more linearly independent operators
$$R_{j+4}(x,\partial_x) = q(x)D(x,\partial_x)^jq(x),\quad\text{for}\ j=0,1,$$
where here linear independence is clear from comparing bandwidths and orders of operators.

In order to complete the construction of a commuting differential operator, we need to find two more non-constant differential operators 
which are linearly independent from the others we have found so far.
To do so, we will use the pair of reciprocal relations
\begin{align*}
\wt{\Psi}(k,x)\cdot q(x)B(x,\partial_x)Q(x,\partial_x)\frac{1}{q(x)} &= 
\frac{1}{p(k)} P(k,\mathscr S_k)A(k,\mathscr S_k)p(k)\cdot\wt{\Psi}(k,x),\\
\wt{\Psi}(k,x)\cdot \frac{1}{q(x)}Q(x,\partial_x)^*B(x,\partial_x)^*q(x) &= 
p(k)A(k,\mathscr S_k)^* P(k,\mathscr S_k)^*\frac{1}{p(k)}\cdot\wt{\Psi}(k,x),
\end{align*}
for all $A(k,\mathscr S_k)\in \mathcal F_k(\Psi)$ and $B(x,\partial_x)\in\mathcal F_x(\Psi)$ with $b_\Psi(A) = B$.
Note that these relations do not preserve formal symmetricity, but if we combine the two we can obtain new formally symmetric operators.
In particular, if $A(k,\mathscr S_k)$ and $B(x,\partial_x)$ are formally skew-symmetric, then for any $j$ we have the formally symmetric difference operator
$$\frac{1}{p(k)}P(k,\mathscr S_k) A(k,\mathscr S_k)^jp(k) + (-1)^j p(k) A(k,\mathscr S_k)^{j*} P(k,\mathscr S_k)^*\frac{1}{p(k)}$$
which is in the Fourier algebra and which the generalized Fourier map sends to the formally symmetric differential operator
$$q(x)B(x,\partial_x)^jQ(x,\partial_x)\frac{1}{q(x)} + (-1)^j \frac{1}{q(x)}Q(x,\partial_x)^*B(x,\partial_x)^{j*}q(x).$$
Taking $B(x,\partial_x) ;= x\partial_x-\partial_xx$, this gives us two more operators
\begin{align*}
R_{j+6}(x,\partial_x) &= q_\pm(x)(x\partial_x-\partial_xx)^jQ_\pm(x,\partial_x)\frac{1}{q_\pm(x)} 
\\
&+ (-1)^j\frac{1}{q_\pm(x)}Q_\pm(x,\partial_x)^*(x\partial_x-\partial_xx)^jq_\pm(x)
\end{align*}
for $j = 0,1$.
The linear independence of these new operators from the previous ones is not immediately clear, but can be verified via computer.
Consequently the symmetric component of the right Fourier algebra of order and co-order $4$, $\mathcal F_{x,sym}^{4,4}(\wt{\Psi})$ is at least $8$ dimensional, containing the constant operator as well as all linear combinations of the seven linearly independent operators $R_1,\dots, R_7$.

The vanishing condition of the discrete and continuous concomitants imposes $6$ conditions, implying the existence of a two dimensional subspace of $\mathcal F_{x,sym}^{4,4}(\wt\Psi)$ satisfying the condition.
Excluding the constants, there must be at least one non-constant operator with this property.
Solving the vanishing condition on the concomitants, we find that the continuous integral operator $T_{\wt\Psi}$ defined by 
Eq. \eqref{eqn:cont integral} will commute with a non-constant differential operator formed by a linear combination of the above seven operators.
Specifically, by solving the associated system via computer, we get that $T_{\wt\Psi}$ commutes with the differential operator of order $4$ defined by
$$R(x,\partial_x) = \sum_{j=1}^7 c_jR_j(x,\partial_x),$$
where
\begin{align*}
c_1 &= t^2,\ c_2 = -2t,\ c_3 = 1,\ c_7 = 0,\\
c_4 &= -\frac{\lambda(\lambda-n)(\lambda+n-1)}{\lambda+a},\\
c_5 &= -\frac{2\lambda(\lambda-n)}{\lambda+a},\\
c_6 &= \frac{\lambda-n}{\beta}\left(a\frac{\beta+1}{\beta-1}-t\right).
\end{align*}

\subsection{A non-scalar example}
\label{7.3}
Noncommutative bispectral Darboux transformations lead to interesting examples of matrix valued discrete-continuous bispectral functions, and thereby examples of matrix valued integral operators commuting with matrix valued differential operators.  As a specific example of this, consider the case when $\Psi(k,x) = \psi(k,x)I$ is the classical discrete-continuous bispectral function of Hermite type and let
$$D(x,\partial_x) := \partial_x^2 - x^2 + 1\quad\text{and}$$
$$L(k,\mathscr S_k) := \sqrt{k/2}\mathscr S_k^{-1} + \sqrt{(k+1)/2}\mathscr S_k$$
be the associated formally symmetric differential and shift bispectral operators.

Fix an $r\times r$ symmetric matrix $A$ and consider the $2r\times 2r$ matrix valued differential operator
$$U(x,\partial_x) := \left(\begin{array}{cc} (\partial_x + x)I & -A\\ -A & (\partial_x-x)I\end{array}\right),$$
where $I$ is the $r\times r$ identity matrix.
The operator $U(x,\partial_x)$ gives rise to a factorization
$$U(x,\partial_x)U(x,\partial_x)^* = \mxx{-D(x,\partial_x)I + A^2}{0}{0}{-D(x,\partial_x)I + A^2 + 2I}.$$
One can feel Dirac lurking in the background.

The right Fourier algebra of $\Psi$ is exactly the algebra of right-acting differential operators with matrix-valued polynomial coefficients.
In fact, the generalized Fourier map satisfies
$$b_\Psi(L) = x,\quad b_\Psi(-2k) = D(x,\partial_x)$$
and also
$$b_\Psi(\sqrt{k/2}\mathscr S_k^{-1} - \sqrt{(k+1)/2}\mathscr S_k) = \partial_x.$$
Therefore the operator $U(x,\partial_x)$ belongs to the right Fourier algebra $\mathcal F_x(\Psi)$.
The associated operator in the left Fourier algebra is
$$P(k,\mathscr S_k) = b_\Psi^{-1}(U(x,\partial_x)) = \left(\begin{array}{cc} \sqrt{2k}\mathscr S_k^{-1}I & -A\\ -A & -\sqrt{2k+2}\mathscr S_kI\end{array}\right)$$
The discrete operator $P$ yeilds the remarkable factorization
$$P(k,\mathscr S_k)^*(F(k)^*)^{-1}F(k)^{-1}P(k,\mathscr S_k) = \left(\begin{array}{cc} I & 0\\0 & I\end{array}\right)$$
for the sequence of $2r\times 2r$ matrices
$$F(k) = \left(\begin{array}{cc} 0 & B_{k}\\B_{k+1} & 0\end{array}\right).$$
Here $B_kB_k^* = 2kI + A^2$ is a fixed Cholesky factorization for each $k$, making
$$F(k)F(k)^* = \left(\begin{array}{cc} 2kI+A^2 & 0\\0 & (2k+2)I+A^2\end{array}\right).$$
This factorization is significant, being of a form appearing in \thref{bisp}.

Thus,
\begin{align*}
\wt{\Psi}(k,x)
  & = F(k)^{-1}\Psi(k,x)\cdot U(x,\partial_x)\\
  & = \left(\begin{array}{cc}-B_{k+1}^{-1}A\psi(k,x) & -B_{k+1}^{-1}\sqrt{2k+2}\psi(k+1,x)\\ B_{k}^{-1}\sqrt{2k}\psi(k-1,x) & -B_{k}^{-1}A\psi(k,x)\end{array}\right).
\end{align*}
defines a self-adjoint noncommutative bispectral Darboux transformation of $\Psi(k,x) = \psi(k,x)I$ satisfying the orthonormality relation
$$\int_{-\infty}^\infty \wt\Psi(j,x)\wt\Psi(k,x)^* dx = \delta_{jk} \left(\begin{array}{cc}I& 0\\0& I\end{array}\right).$$

As a consequence, \thref{gral} implies that the continuous integral operator $T_{\wt\Psi}$ defined by Eq. \eqref{eqn:cont integral}
$$T_{\wt\Psi}: F(y)\mapsto \int_{-\infty}^{t} F(x)K(x,y) dx,\ \ K(x,y) = \sum_{k=0}^{n} \wt\Psi(k, x)^*\wt\Psi(k,y)$$
which acts on matrix-valued functions will commute with a matrix valued differential operator of positive order.

The bispectral function $\wt\Psi(x,y)$ satisfies a three-term recursion relation, allowing us to obtain an explicit expression for $K(x,y)$.
In particular,
$$x\wt\Psi(k,x) = H_1(k)\wt\Psi(k+1,x) + H_0(k)\wt\Psi(k,x) + H_1(k-1)^*\wt\Psi(k-1,x)$$
for the matrices
\begin{align*}
H_1(k) &= \sqrt{(k+1)/2}\cdot F(k)^{-1}\left(\begin{array}{cc}A^2+2kI&0\\0& A^2+(2k+4)I\end{array}\right)(F(k+1)^*)^{-1}\\
       &= \sqrt{(k+1)/2}\left(\begin{array}{cc}B_{k+1}^{-1}B_{k+2}& 0 \\ 0 & B_k^*(B_{k+1}^*)^{-1}\end{array}\right)\\
H_0(k) &= F(k)^{-1}\left(\begin{array}{cc}0&A\\A& 0\end{array}\right)(F(k)^*)^{-1} = \left(\begin{array}{cc}0 & B_{k+1}^{-1}A(B_k^*)^{-1}\\ B_k^{-1}A(B_{k+1}^*)^{-1} & 0\end{array}\right).
\end{align*}
Therefore the kernel $K(x,y)$ may be expressed as
$$K(x,y) = \frac{1}{x-y}\left(\wt\Psi(n+1,x)^*H_1(n)\wt\Psi(n,y) - \wt\Psi(n,x)^*H_1(n)\wt\Psi(n+1,y)\right).$$

To determine the commuting operator, we must explore the Fourier algebras and the generalized bispectral map of $\wt\Psi$.
The right Fourier algebra of $\wt\Psi$ contains every matrix-valued differential operator with polynomial coefficients.  In fact, if $R(x,\partial_x)$ is such an operator then its preimage under the generalized Fourier map is
$$b_{\wt\Psi}^{-1}(R(x,\partial_x)) = F(k)^{-1}P(k,\mathscr S_k)b_\Psi^{-1}(R(x,\partial_x))P(k,\mathscr S_k)^*(F(k)^*)^{-1}.$$
On the other hand for the left Fourier algebra, if $M(k,\mathscr S_k)\in \mathcal F_k(\Psi)$ then
$$b_{\wt\Psi}(F(k)^*M(k,\mathscr S_k)F(k)) = U(x,\partial_x)^*b_\Psi(M(k,\mathscr S_k)) U(x,\partial_x).$$

We can get an estimate of the dimension of $\mathcal F^{4,4}_{sym,x}(\wt\Psi)$ by first obtaining appropriate estimates of the vector spaces of scalar operators $\mathcal F_{x}^{2,4}(\psi)$ and $\mathcal F_{x}^{1,4}(\psi)$.
Direct computation shows that the symmetric components of these spaces are
$$\mathcal F_{x,sym}^{1,4}(\psi) = \Span\{1,x,x^2\}$$
and
$$\mathcal F_{x,sym}^{2,4}(\psi)  = \Span\{1,x,x^2,D(x,\partial_x),\{x,D(x,\partial_x)\},\{x^2,D(x,\partial_x)\}\}.$$
The skew-symmetric component of each space is the same, given by
$$\mathcal F_{x,skew}^{2,4}(\psi) = \mathcal F_{x,skew}^{1,4}(\psi) = \Span\{[x,D(x,\partial_x)],[x^2,D(x,\partial_x)]\}.$$
Here $\{\cdot,\cdot\}$ denotes the anticommutator and $[\cdot,\cdot]$ denotes the commutator.

Therefore the symmetric component $\mathcal F_{x,sym}^{4,4}(\wt\Psi)$ of the right Fourier algebra of $\wt\Psi(k,x)$ contains any operator of the form
\begin{align*}
R(x,\partial_x)
  & = U(x,\partial_x)^*\left(\sum_{j=0}^2\sum_{k=0}^1A_{j,k}\{x^j,D(x,\partial_x)^k\} + \sum_{j=1}^2B_j[x^j,D(x,\partial_x)]\right)U(x,\partial_x)\\
  & + V_0 + V_1x + V_2x^2 + W_1[x,D(x,\partial_x)] + W_2[x^2,D(x,\partial_x)]
\end{align*}
for symmetric matrices $A_j,V_j$ and skew-symmetric matrices $B_j,W_j$, where $V_2$ and $W_2$ satisfy a certain set of $N=2r$ conditions to force $b^{-1}(R)$ to have bandwidth $\leq 4$.
In particular, this gives us the dimension estimate
$$\dim\mathcal F_{x,sym}^{4,4}(\wt\Psi) \geq 9\frac{N(N+1)}{2} + 4\frac{(N-1)N}{2} = 13N^2/2 + 3N/2.$$
The vanishing condition on the left and right concomitants imposes $6N^2$ conditions.
Consequently, there exist at least $N^2/2 + 3N/2$ operators commuting with the integral operator $T_{\wt\Psi}$ defined above.
By \thref{bisymmetric existence}, this implies the existence of a differential operator of order no greater than four commuting with $T_{\wt\Psi}$.
At most $N(N+1)/2$ of these operators are constant, so this implies that the vector space $\mathcal F_{x,sym}^{4,4}(\wt\Psi)$ contains at least $N$ non-constant operators commuting with $T_{\wt\Psi}$.
An exact expression for a commuting operator can be found by solving for the vanishing condition of the concomitants.
Solving for this condition, we find that in fact the operator $T_{\wt\Psi}$ commutes with the \emph{second-order differential operator}
\begin{align*}
R(x,\partial_x)
  & = U(x,\partial_x)^*(t-x)U(x,\partial_x) + x(A^2+(2n+2)I)I_N\\
  & = \partial_x (x-t)\partial_xI_N + (tA^2 + (2n+2)x+tx^2-x^3)I_N + \left(\begin{array}{cc}t-2x&A\\A&-t+2x\end{array}\right),
\end{align*}
where $I_N = \left(\begin{array}{cc}I& 0\\0&I\end{array}\right)$ is the $N\times N$ identity matrix.

\subsection{An example not coming from orthogonal polynomials}
From the examples depicted so far, it is tempting to conclude that all examples are obtainable 
from bispectral Darboux transformations of the classical discrete-continuous bispectral functions.
As the following example shows, this is not the case.
In particular there are examples of a discrete-continuous bispectral functions arising form ``pure soliton" solutions of the KdV equation.

Let $N>0$ be a positive integer and choose positive constanst $k_1,\dots, k_N$ and $c_1,\dots, c_N$ and define matrices
\begin{align*}
M_{ij}(x) &:= \delta_{ij} + \frac{c_j}{k_i+k_j}e^{(k_i+k_j)x},
\\
\wt M_{ij}(k,x) &:= \delta_{ij} + \frac{c_j}{k_i+k_j}\frac{k-k_j}{k+k_j}e^{(k_i+k_j)x}.
\end{align*}
Then the functions $\tau(x) := \det M(x)$ and $\tau(k,x) := \det\wt M(k,x)$ define the wave function
$$\Psi(k,x) := e^{kx}\frac{\tau(k,x)}{\tau(x)}$$
satisfying the Schr\"{o}dinger equation
$$\Psi(k,x)\cdot D(x,\partial_x) = k^2\Psi(k,x),\quad\text{for}\ D(x,\partial_x) := \partial_x^2 + V(x),$$
where the potential $V(x)$ is given by $V(x) := 2(\log \tau(x))''$.

Now if we specifically choose $k_j := j$ for all $j=1,\dots, N$ and
$$c_j := 2j(-1)^{j-1}\prod_{i\neq j}\frac{i+j}{i-j},$$
then $V(k) = N(N+1)\sech^2(k)$ and a solution of the above Schr\"{o}dinger equation can be expressed in terms of the hypergeometric function
$$\Psi(k,x) = \mu(k)e^{kx}\ _2F_1\left(-N,N+1,1+k; \frac{e^x}{e^x-e^{-x}}\right),$$
where here for sake of normalization we define $\mu(k)$ recursively by 
\[
\mu(1) := 1 \quad \mbox{and} \quad \mu(k+1)/\mu(k) := A(k) \; \; \mbox{for} \; \; A(k) := \left(\frac{(N-k)(k+N+1)}{k(k+1)}\right)^{1/2}.
\]
This function satisfies the difference equation 
\[
L(k,\mathscr S_k)\cdot\Psi(k,x) = 2\sinh(x)\Psi(k,x)
\]
for
$$L(k,\mathscr S_k) = A(k)\mathscr S_k + A(k-1)\mathscr S_k^{-1}.$$

For $k=1,2,\dots, N$, the function $\Psi(k,x)$ is square-integrable on $\mathbb R$, and by the choice of $\mu(k)$ above has norm $1$.
The integral operator
$$T_\Psi(F)(z) = \int_t^\infty F(x)K(x,z) dx$$
with kernel
$$K(x,z) = \sum_{j=p}^N \Psi(j,x)\Psi(j,z)$$
can be viewed as an integral operator on $L^2(\mathbb R)$.

To obtain a differential operator commuting with this integral operator, we wish to get a better sense of the operators in the symmetric part of the Fourier algebras $\mathcal F_{k,sym}(\Psi)$ and $\mathcal F_{z,sym}(\Psi)$.
We showed above that 
\[
L(k,\mathscr S_k)\cdot \Psi(k,x) = \Psi(k,x)2\sinh(x) \quad \mbox{and} \quad \Psi(k,x)\cdot D(x,\partial_x) = k^2\Psi(k,x),
\]
and therefore,
\[
k^2,L(k,\partial_k)\in \mathcal F_{k,sym}(\Psi) \quad \mbox{and} \quad 2\sinh(x),D(x,\partial_x)\in \mathcal F_{x,sym}(\Psi).
\]
More generally, the anticommutator of symmetric operators is also symmetric, so
\begin{align*}
\{k^2,L\}
  & = k^2L(k,\mathscr S_k) + L(k,\mathscr S_k)k^2\\
  & = A(k)(2k^2+2k+1)\mathscr S_k + A(k-1)(2k^2-2k+1)\mathscr S_k^{-1}
\end{align*}
is an element of $\mathcal F_{k,sym}(\Psi)$ and 
\begin{align*}
\{2\sinh(x),D(x,\partial_x)\}
&= 2\sinh(x)D(x,\partial_x) + 2D(x,\partial_x)\sinh(x)\\
&= 4\partial_x(\sinh(x))\partial_x + 2\sinh(x) + 4\sinh(x)N(N+1)\sech^2(x)
\end{align*}
is an element of $\mathcal F_{x,sym}(\Psi)$.

Since the generalized Fourier map is an algebra isomorphism, it must map anticommutators to anticommutators, and thus,
$$b_\Psi: \left.\begin{array}{ccc}
L(k,\mathscr S_k)&\mapsto &2\sinh(x),\\
k^2 &\mapsto & D(x,\partial_x),\\
\{k^2,L\} & \mapsto & \{2\sinh(x),D(x,\partial_x)\}.
\end{array}\right.$$
More explicitly, $b_\Psi$ satisfies
\begin{align*}
A(k)\mathscr S_k + A(k-1)\mathscr S_k^{-1}& \mapsto 2\sinh(x),\\
k^2 & \mapsto \partial_x^2 + N(N+1)\sech^2(k),\\
 \left(\begin{array}{l}A(k)(2k^2+2k+1)\mathscr S_k\\\quad+ A(k-1)(2k^2-2k+1)\mathscr S_k^{-1}\end{array}\right)
& \mapsto \left(
\begin{array}{l}4\partial_x(\sinh(x))\partial_x + 2\sinh(x) \\\quad+ 4\sinh(x)N(N+1)\sech^2(x)\end{array}\right),
\end{align*}
It follows that for $*$ the anti-involution, $\mathcal F_{k,sym}^{2,2}(\Psi)$ has dimension $8$ as a real vector space (including the constants).

By \thref{bisymmetric existence} combined with \thref{basic}, there must exist a nonconstant differential operator in $\mathcal F_{y,sym}^{2,2}(\Psi)$ commuting with the integral operator $T_\Psi$.
To find it, we can take a general linear combination of the non-constant operators in $\mathcal F_{k,sym}^{2,2}(\Psi)$:
$$ P (k,\mathscr S_k) := A(k)(2k^2+2k+1+\beta)\mathscr S_k +\gamma k^2+ A(k-1)(2k^2-2k+1+\beta)\mathscr S_k^{-1}.$$
By linearity, the image $R(x,\partial_x) := b_\Psi(P)$ satisfies
\begin{align*}
R(x,\partial_x) & = \partial_x(4\sinh(x) + \gamma)\partial_x+ 2\sinh(x)(\beta+1)\\
                &  +4\sinh(x)N(N+1)\sech^2(x) + \gamma N(N+1)\sech^2(x).
\end{align*}
The general bi-self-adjoint operator we are looking for is exactly the difference operator $L$ whose discrete concomitant \eqref{discrete concomitant} vanishes at $p-1$, and whose continuous concomitant vanishes at $t$.
The first condition implies that $\beta = -2(p-1)^2-2(p-1)-1$, while the second condition gives $\gamma = -4\sinh(t)$.
Inserting these values and rescaling the expression by a factor of $4$, we find that the differential operator
\begin{align*}
R(x,\partial_x) &= \partial_x(\sinh(x)-\sinh(t))\partial_x - p(p-1)\sinh(x)\\
&+ N(N+1)(\sinh(x)-\sinh(t))\sech^2(x)
\end{align*}
commutes with the integral operator $T_\Psi$.


\begin{thebibliography}{xx}
\bibitem{BHY1} B.~Bakalov, E.~Horozov, and M.~Yakimov, {\em{General methods for constructing bispectral operators}}, 
Phys. Lett. A {\bf{222}} (1996), 59--66.

\bibitem{BHY2} B.~Bakalov, E.~Horozov, and M.~Yakimov, {\em{Bispectral algebras of commuting ordinary differential operators}}, 
Comm. Math. Phys. {\bf{190}} (1997), 331--373.

\bibitem{bochner1929sturm}
S. Bochner, {\em{{\"U}ber {S}turm-{L}iouvillesche polynomsysteme}}, 
Math. Z. {\bf{29}} (1929), 730--736.

\bibitem{CGYZ1} W. R. Casper, F. A. Gr\"unbaum, M. Yakimov, and I. Zurri\'an, {\em{Reflective prolate-spheroidal operators and the KP/KdV equations}}, 
Proc. Natl. Acad. Sci. USA {\bf{116}} (2019), no. 37, 18310--18315.




\bibitem{CGYZ3} W. R. Casper, F. A. Gr\"unbaum, M. Yakimov, and I. Zurri\'an, {\em{Algebras of commuting differential operators for kernels of Airy type}}, 
Toeplitz Operators and Random Matrices: In Memory of Harold Widom. Springer International Publishing, 2022. 229--256.

\bibitem{CGYZ2} W. R. Casper, F. A. Gr\"unbaum, M. Yakimov, and I. Zurri\'an, {\em{Reflective prolate-spheroidal operators and the adelic Grassmannian}},
Comm. Pure Appl. Math., {\bf{76}} (2023), no. 12, 3769--3810. 

\bibitem{CY1} W. R. Casper and  M. Yakimov, {\em{ Integral operators, bispectrality and growth of Fourier algebras}}, J. Reine Angew. Math. {\bf 766} (2020), 151--194.

\bibitem{CY2} W. R. Casper and M. Yakimov, {\em{The matrix {B}ochner problem}}, Amer. J. Math. {\bf 144} (2022), no. 4, 1009--1065.

\bibitem{coddington} E. A. Coddington and N. Levinson, {\em{Theory of ordinary differential equations}}, Tata McGraw-Hill Education, 1955.

\bibitem{C} A. Connes, {\em{Trace formula in noncommutative geometry and the zeros of the Riemann zeta function}},  
Selecta Math. (N.S.) {\bf{5}} (1999), no. 1, 29–106.

\bibitem{CC} A. Connes and C. Consani, {\em{Spectral triples and $\zeta$-cycles}}, arXiv:2106.01715.

\bibitem{CM} A. Connes and H. Moscovici, {\em{The UV prolate spectrum matches the zeros of zeta}}, 
Proc. Natl. Acad. Sci. U.S.A. {\bf{119}}, e2123174119 (2022).

\bibitem{DG} J.~J. Duistermaat and F.~A. Gr{\"{u}}nbaum, {\em{Differential equations in the spectral parameter}}, 
Comm. Math. Phys. {\bf{103}} (1986), 177--240.

\bibitem{DuG} A. J. Dur\'an and F. A. Gr\"unbaum, {\em{Orthogonal matrix polynomials satisfying second-order differential equations}}, 
Int. Math. Res. Not. {\bf{2004}}, no. 10, 461--484.

\bibitem{GHY} J. Geiger, E. Horozov, and M. Yakimov, {\em{Noncommutative bispectral Darboux transformations}}, Trans. Amer. Math. Soc. {\bf{369}} (2017), 
5889--5919.

\bibitem{grunbaum1}
F. A. Gr\"{u}nbaum, {\em{{A} new property of the reproducing kernels for classical orthogonal polynomials}},
J. Math. Anal. Appl. {\bf{95}} (1983), no. 2, 491--500.


\bibitem{GH} F. A. Gr\"unbaum and L. Haine, {\em{Some functions that generalize the Askey--Wilson polynomials}}, 
Comm. Math. Phys. {\bf{184}} (1997), no. 1, 173--202.

\bibitem{GHH} F. A. Gr\"unbaum, L. Haine, and E. Horozov, {\em{On the Krall--Hermite and the Krall--Bessel polynomials}}, 
Internat. Math. Res. Notices {\bf{1997}}, no. 19, 953--966.

\bibitem{GPT} F. A. Gr\"unbaum, I. Pacharoni, and J. Tirao, 
{\em{Matrix valued spherical functions associated to the complex projective plane}}, 
J. Funct. Anal. {\bf{188}} (2002), 350--441. 

\bibitem{GPZ} F. A. Gr\"unbaum, I. Pacharoni, and I. Zurri\'an, {\em{Bispectrality and time-band limiting: matrix-valued polynomials}}, 
Int. Math. Res. Not. IMRN {\bf{2020}}, no. 13, 4016--4036.

\bibitem{HL} L. Haine and P. Iliev, {\em{Commutative rings of difference operators and an adelic flag manifold}}, 
Internat. Math. Res. Notices {\bf{2000}}, no. 6, 281--323.

\bibitem{KR} A. Kasman and M. Rothstein, {\em{Bispectral Darboux transformations: the generalized Airy case}},
Phys. D {\bf{102}} (1997), 159--176. 

\bibitem{K} M. G. Kre\u{\i}n, {\em{Infinite $J$-matrices and a matrix-moment problem}}, Doklady Akad. Nauk SSSR (N.S.) {\bf{69}} (1949), 125--128.

\bibitem{LP} H.~J. Landau and H.~O. Pollak, {\em{Prolate spheroidal wave functions, {F}ourier analysis and
uncertainty. {II}}}, Bell System Tech. J. {\bf{40}} (1961), 65--84.

\bibitem{Mehta} M.~L. Mehta, {\em Random matrices}, 3rd ed, Elsevier Publ., 2004.

\bibitem{reach2} M. Reach, {\em{Generating difference equations with the {D}arboux transformation}}, Comm. Math. Phys. {\bf 119} (1988), no 3, 385--402.

\bibitem{reach3} M. Reach, {\em{Difference equations for {$N$}-soliton solutions to {K}d{V}}}, Phys. Lett. A {\bf 129} (1988), no 2, 101--105.

\bibitem{S} D.~Slepian, {\em{Prolate spheroidal wave functions, {F}ourier analysis and uncertainity. {IV}. {E}xtensions to many dimensions; generalized prolate
spheroidal functions}}, Bell System Tech. J. {\bf{43}} (1964), 3009--3057.

\bibitem{SP} D.~Slepian and H.~O. Pollak, {\em{Prolate spheroidal wave functions, {F}ourier analysis and
uncertainty. {I}}},  Bell System Tech. J. {\bf{40}} (1961), 43--63.

\bibitem{TW1} C.~A. Tracy and H.~Widom, {\em{Fredholm determinants, differential equations and matrix models}}, 
Comm. Math. Phys. {\bf{163}} (1994), 33--72.

\bibitem{TW2} C.~A. Tracy and H.~Widom, {\em{Level-spacing distributions and the {A}iry kernel}}, 
Comm. Math. Phys. {\bf{159}} (1994), 151--174.

\bibitem{Wilson} G. Wilson, {\em{Bispectral commutative ordinary differential operators}}, J. Reine Angew. Math. {\bf{442}} (1993), 177--204.


\bibitem{Z} J. P. Zubelli, {\em{Differential equations in the spectral parameter for matrix differential operators}}, 
Phys. D {\bf{43}} (1990), no. 2-3, 269--287. 



\end{thebibliography}
\end{document}